\documentclass [12pt,leqno]{article}
\usepackage[leqno]{amsmath}

\usepackage{enumerate,amsthm,amsfonts,amsmath,amssymb,graphics, pst-poly,pstricks}%

\divide\oddsidemargin
by 2  

\makeatletter
\newcommand{\address}[2][]{%
  \ifx\@add@ress\@undefined\gdef\@add@ress{\par\par\bigskip}\AtEndDocument{\@add@ress}\fi%
  \g@addto@macro\@add@ress{\bigskip\noindent{\small\scshape%
      \ifx#1\empty\else{\bfseries Address of #1:}\ \fi#2}\par\par}}
\newcommand{\email}[1]{\par\par\noindent{\normalfont{\itshape E-mail:\/} \texttt{#1}}\par\par}
\renewenvironment{abstract}{\small\quotation\noindent
  {\bfseries \abstractname}}{\endquotation \par}
\newcommand{\footnotetextplain}[1]{\begingroup\def\@thefnmark{}%
  \long\def\@makefntext##1{\parindent 0pt\noindent ##1}\@footnotetext{#1}
  \endgroup}
\newcommand{\AMSsubjclass}[2]{\footnotetextplain{2000
   \emph{Mathematics Subject Classification:} Primary #1, Secondary #2.}}

\newcommand{\keywords}[1]{\footnotetextplain{\emph{Key words and phrases:} #1.}}
\makeatletter \xdef\qedbuit{\qed}

\newcommand{\TeoremaAmbFinalMarcat}[1]{%
  \expandafter\gdef\csname end#1\endcsname{\qedbuit\@endtheorem}}
\newtheorem{theo}{Theorem}[section]

\theoremstyle{definition}
 \TeoremaAmbFinalMarcat{defi}
 \TeoremaAmbFinalMarcat{ex}
\newtheorem{rem}[theo]{Remark} \TeoremaAmbFinalMarcat{rem}
 \TeoremaAmbFinalMarcat{conv}

\newenvironment{proclama}[1]{\par\vspace{\topsep}\noindent{\bf #1}
 \begin{em}}{\end{em}\par\vspace{\topsep}}

\newcommand{\start}[2]{\begin{#1}\label{#2}}

\newcommand{\theoc}[1]{Theorem~\ref{#1}}
\newcommand{\propc}[1]{Proposition~\ref{#1}}
\newcommand{\coryc}[1]{Corollary~\ref{#1}}
\newcommand{\lemc}[1]{Lemma~\ref{#1}}

\newcommand{\remc}[1]{Remark~\ref{#1}}
\newcommand{\exc}[1]{Example~\ref{#1}}
\newcommand{\figc}[1]{Figure~\ref{#1}}

\newcommand{\defc}[1]{Definition~\ref{#1}}
\newcommand{\refc}[1]{~\ref{#1}}

\def\@enum@{\list{\csname label\@enumctr\endcsname}%
           {\usecounter{\@enumctr}\def\makelabel##1{\hss\llap{##1}}
           \itemsep=2pt\parsep=0pt\topsep=3pt plus 1pt minus 1 pt}}
\newenvironment{alphlist}{\enumerate[(a)]}{\endenumerate}
\newenvironment{romlist}{\enumerate[(i)]}{\endenumerate}
\newenvironment{numlist}{\enumerate[(1)]}{\endenumerate}

\newcommand{\sequence}{%
{\rput{18}(1.4,1.4){\PstPolygon[unit=1.3,PolyNbSides=10] }
\rput{20}(3.9,1.4){\PstPolygon[unit=1.3,PolyNbSides=10] }
\rput{20}(6.3,1.4){\PstPolygon[unit=1.3,PolyNbSides=10]}
\rput{20}(8.8,1.4){\PstPolygon[unit=1.3,PolyNbSides=10]}
\rput{20}(11.3,1.4){\PstPolygon[unit=1.3,PolyNbSides=10]} } }

\newcommand{\TorusWithOneHole}{
\pscurve[](9,3)(5,3)(5,0.6)(9,0.6)
\psccurve(9,3)(8.8,2.6)(8.8,1)(9,0.6)(9.2,1)(9.2,2.6)
\pscurve(6,2)(7,1.8)(8,2)
\pscurve(7.8,1.95)(7,2.1)(6.7,2.1)(6.2,1.96) }

\newcommand{\EdgesForTorusWithOneHole}{
\pscurve[linecolor=lightgray](8.8,2.6)(6,3)(5.2,2)(6,1)(8.8,1.2)
\pscurve[linecolor=lightgray](8.8,2.4)(8.2,2.3 )(7.85,2) \pscurve[linestyle=dotted,
dotsep=1.8pt](7.8,2)(8.1,1.7)(8.7,1.6)
\psecurve[linecolor=lightgray](8.1,1.7)(8.7,1.6)(9.3,1.6)(9.5,1.6)}
\newcommand{\CurvesForTorusWithOneHole}{
\psset{linecolor=darkgray}
\psecurve[showpoints=false]{->}(7,0)(7.4,3.3)(8.4,1.5)(5.4,1.9)(8.1,2.2
)(7.9,1.83)(7.6,1.83)(7.3,1.78)(7.3,1.78)
\psecurve[linestyle=dashed]{->}(8,1)(7.3,1.78)(6.8,0.2)(6.2,-1)
\psecurve{->}(7,4)(6.8,0.2)(6,1)(5.8,2)(6.2,2.2)(6.69,2.1)(17,2)
\psecurve[linestyle=dashed]{->}(4,4)(6.6,2.1)(7.4,3.3)(7,0)
 }

\def\flecha#1#2#3{\mbox{${#1}\stackrel{#2}{\longrightarrow} {#3}$}}

\def\map#1#2#3{\mbox{${#1}\colon {#2} \longrightarrow {#3}$}}
\def\Smap#1#2{\mbox{${#1}\colon{#2} \longrightarrow {#2}$}}
\def\conj#1{\mbox{$\{#1\}$}}
\def\ai{\mathsf{I}}
\def\ap{\mathsf{P}}
\def\aq{\mathsf{Q}}

\def\VVV{\mathsf{V}}
\def\YYY{\mathsf{Y}}
\def\XXX{\mathsf{X}}
\def\WWW{\mathsf{W}}
\def\ZZZ{\mathsf{Z}}
\def\CC{\mathcal{C}}
\def\XX{\mathcal{X}}
\def\UU{\mathcal{U}}
\def\VV{\mathcal{V}}
\def\WW{\mathcal{W}}
\def\ZZ{\mathcal{Z}}
\def\OO{\mathcal{O}}
\def\ov#1{\overline{#1}}
\def\ova{\overline{a}}
\def\ovx{\overline{x}}
\def\ovy{\overline{y}}
\def\bb{\mathbf{LP}}
\def\id{\mathop\mathrm{Id}}
\def\V{\mathbb{V}}
\def\len{\mathsf{l}}
\def\sign{\mathrm{sign}}
\def\si{\mathop\mathrm{o}}
\def\class{\mathrm{class}}
\def\im{\mathrm{Im}}

\def\cob{\delta}
\def\cib{\Delta}
\def\br{\gamma}
\def\bra#1{[#1]}
\def\bri#1{\langle #1 \rangle}
\def\om{\omega}
\def\sss{\mathrm{s}}

\def\Po{P_\mathcal{O}}
\def\So{\Sigma_\mathcal{O}}
\newcommand{\an}{\ensuremath{\mathbb{A}}_n}
\newcommand{\atwo}{\ensuremath{\mathbb{A}}_2}

\def\ox{\overline{x}}

\def\cc{\mathsf{c}}
\def\DD{\mathsf{D}}
\def\W{\mathbb{W}}
\def\S{\mathsf{S}}
\newcommand{\N}{\ensuremath{\mathbb{N}}}
\newcommand{\Z}{\ensuremath{\mathbb{Z}}}
\title{Combinatorial Lie bialgebras of curves on surfaces}
\author{Moira Chas}
\address{Institute for Mathematical Sciences
Stony Brook University Stony Brook, NY 
11794-3660\email{moira@math.sunysb.edu}} 
\date{\empty} 

\begin{document}

\maketitle

\begin{abstract} Goldman \cite{Gol} and Turaev \cite{T} found a Lie bialgebra structure
on the vector space generated by non-trivial free homotopy classes of curves on a
surface. When the surface has non-empty boundary, this vector space has a basis of cyclic
reduced words in the generators of the fundamental group and their inverses. We give a
combinatorial algorithm to compute this Lie bialgebra on this vector space of cyclic
words. Using this presentation, we prove a variant of Goldman's result relating
the bracket to disjointness of curve representatives when one of the classes is simple.
We exhibit some examples we found by programming the algorithm which answer negatively
Turaev's question about the characterization of simple curves in terms of the cobracket.
Further computations suggest an alternative characterization of simple curves in terms of
the bracket of a curve and its inverse. Turaev's question is still open in genus zero.
\end{abstract}

\AMSsubjclass{57M99}{17B62}

\keywords{surfaces, conjugacy classes, Lie bialgebras}

\section{Introduction}

A Lie bialgebra structure on vector space $W$ consists in two 
linear operations, a bracket from $W\otimes W$ to $W$ and a 
cobracket, from $W$ to $W\otimes W$, satisfying certain 
identities (see Appendix \ref{definition of the Lie 
bialgebra}). Goldman \cite{Gol} and Turaev \cite{T} found, in 
stages, a Lie bialgebra structure on  the vector space 
generated by all non-trivial free homotopy classes of curves 
on an orientable surface. The desire to understand better the 
beautiful structure of Goldman and Turaev and to answer some 
of the questions posed by them motivated this work. The Lie 
algebra of Goldman, as well as the Goldman-Turaev Lie 
bialgebra, can be generalized via "String topology" to 
manifolds of all dimensions, see \cite{CS1} and \cite{CS2}. 

Here,
\begin{itemize}
\item we give explicit presentations of the Goldman-Turaev Lie bialgebra of curves on a
surface with boundary, (there is one for each surface symbol, see Section \ref{words}).
\end{itemize}
These presentations define purely combinatorial Lie bialgebra structures on the vector
space of reduced cyclic words on certain alphabets and, therefore, give  algorithms to
compute algebraically the bracket and cobracket. These algorithms can be programmed and
so we did, finding examples which answer certain questions about the Goldman-Turaev Lie
bialgebra we describe now.

Goldman \cite{Gol} showed that if the bracket of the two free 
homotopy classes is zero, and one of them has a simple 
representative, then these classes have disjoint 
representatives. His proof uses Kerckhoff's convexity 
property of Teichmuller space \cite{Ke}. We extend this 
result by showing that 
\begin{itemize}
\item the number of terms of the bracket of two classes, one of them simple and non-homologous to zero, equals the
minimum number of intersection points of these classes 
(\theoc{no cancellation}).
\end{itemize}                                                          

\begin{proclama}{Theorem \ref{no cancellation}} Let $\VV$ and $\WW$ be cyclic reduced words and  such
that $\VV$ has a simple representative which is non-homologous to zero. Then there exists two  representatives $\alpha$
and $\beta$ of $\VV$ and $\WW$ respectively such that the bracket of $\VV$ and $\WW$
computed using the intersection points of  $\alpha$ and $\beta$ does not have
cancellation. In other words, the number of terms (counted with multiplicity) of
$\bri{\VV, \WW}$ equals the minimal number of intersection points of representatives of
$\VV$ and $\WW$.

\end{proclama}

On the last page of \cite{Gol}, Goldman asked whether it 
would be possible to replace Kerckhoff convexity by a 
topological argument, 
\begin{itemize}
\item the proof here of the variant of Goldman's result is essentially topological.
\end{itemize}
The hypothesis that one of the classes is simple cannot be omitted. By running our
program, we found out that
\begin{itemize}
\item there exist pairs of distinct classes, which are not multiples of simple curves,
have bracket zero and do not have disjoint representatives (\exc{and more}).
\end{itemize}
Goldman recently found these examples independently (email communication),
see also the
last page of \cite{Gol}.

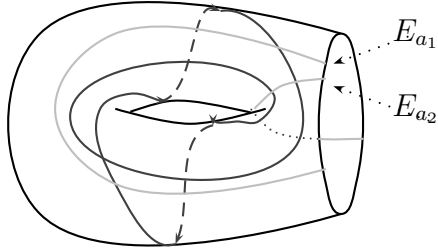
\begin{figure}[ht]\nonumber
\begin{pspicture}(15.5,4)
\TorusWithOneHole \CurvesForTorusWithOneHole 
\EdgesForTorusWithOneHole 
\rput(10,3){$E_{a_1}$}\rput(10,2){$E_{a_2}$} 

\psline[linestyle=dotted, linecolor=black]{->}(10,3)(8.9,2.6) 

\psline[linestyle=dotted, 
linecolor=black]{->}(10.28,1.8)(8.9,2.3) 
\end{pspicture}\nonumber
\caption{A representative of $\cc(a_1a_1 a_2 a_2)$ in the 
punctured torus}\label{counterexample1} 
\end{figure}

In \cite{T}, Turaev formulated a statement dual to that of Goldman's, that is, if the
cobracket of a class is zero, then the class is a multiple of a simple curve and asked
whether this is true.  Again, by the aid of the computer, we found that
\begin{itemize}
\item in every surface of negative Euler characteristic and positive genus
there exists classes with cobracket zero which are not 
multiples of simple curves (\exc{counter}, 
Figure~\ref{counterexample1}). 
\end{itemize}

\begin{itemize}
\item Turaev's conjectural characterization of simple 
curves is still possible for genus zero surfaces.

\item a possible replacement for Turaev's condition in surfaces of all genus
is that multiples of simple curves are characterized by the 
vanishing of the bracket of a class with its inverse. 
Moreover, the output of our program suggests a stronger 
statement: the number of terms of the bracket of a primitive 
class with its inverse is twice the minimal number of 
self-intersection points of the class.

In fact we have quantitative results about the last two 
questions. 
\end{itemize}

\begin{proclama}{Theorem~\ref{new}}(1) On the 
sphere with three punctures all the cyclic words with at most 
sixteen letters, except the multiples of the three peripheral 
curves, have non-zero cobracket. 

(2) On the torus with two punctures all the cyclic words 
$\alpha$ with at most fifteen letters have the property  that 
twice the minimal number of self-intersection points equals 
the number of terms of the bracket 
$[\alpha,\overline{\alpha}]$ in the natural basis. 

\end{proclama}

Some examples we have computed suggest that even a more 
general result may hold: 

\begin{proclama}{Question.} Let $n$ and $m$ be two different non-zero integers and let
$\VV$ be a primitive reduced cyclic word. Is the number of 
terms of the bracket of $\VV^n$ with $\VV^m$ equal to 
$2|m.n|$ multiplied by the minimal number of 
self-intersection points of $\VV$?\end{proclama}

Now we wonder about the possible implications of the computer 
program for counterexamples in three-manifold theory using 
complicated collections of disjoint simple curves to generate 
Heegard decompositions. 

A new version of the programs to compute bracket and cobracket, 
and to  compute  intersection numbers of curves on orientable surfaces with
boundary can be found in the author's homepage: http://www.math.sunysb.edu/~moira/

The rest of the paper is organized as follows. In Section~\ref{words}, we define $\V$,
the vector space of reduced cyclic words on certain alphabets, fix a surface symbol, and
associate to each reduced cyclic word $\VV$, a certain subset $\bb_1(\VV)$ of pairs of
subwords which will play a key role later. Analogously, to each pair of reduced cyclic
words $\VV$ and $\WW$ we associate a subset $\bb_2(\VV,\WW)$ of pairs of the form
(subword of $\VV$, subword of $\WW$). Using $\bb_1(\VV)$ we define a linear map from $\V$
to $\V \otimes \V$ and using
 $\bb_2(\VV,\WW)$ we define a linear map from $\V \otimes \V$ to $\V$.

Free homotopy classes of curves on a surface with boundary, 
the latter described by a surface symbol in an appropriate 
alphabet, are in one-to-one correspondence with cyclic 
reduced words in that alphabet, that is, with the basis of 
$\V$. Since the Goldman-Turaev Lie bialgebra is defined using 
intersection points of representing geometric curves, in 
Section~\ref{Intersection} we study the relation between 
cyclic words and intersection points of certain 
representatives. More precisely, we show that each primitive 
cyclic word $\VV$ in  $\V$ has a representative such that its 
self-intersection points are in one-to-one correspondence 
with $\bb_1(\VV)$ quotiented by an involution 
(\theoc{bigons}). Furthermore, for each pair of cyclic words 
$\VV, \WW$ there exists a pair of representatives such that 
the intersection points of $\VV$ and $\WW$ are in one-to-one 
correspondence with elements of $\bb_2(\VV,\WW)$ 
(\theoc{bigons for a pair}). 

Using the correspondence between free homotopy classes and cyclic reduced words, the
Goldman-Turaev Lie bialgebra operations become defined on $\V$. In Section~\ref{GT},
using the representatives of Section~\ref{Intersection}, we prove that two linear maps we
define in Section~\ref{words} are the Goldman-Turaev bracket and cobracket
(\propc{isomorphism}).

In Section~\ref{applications} we give a topological proof of the variant of the result of Goldman about
simple curves and disjointness (\coryc{cory}) and a generalization of this variant (\theoc{no
cancellation}), and we exhibit some examples that shows that a dual statement to that of
Goldman asked by Turaev does not hold (\exc{counter}). We conclude the section by stating
open problems related to Turaev's characterization and its possible replacement.

In  Appendix~\ref{definition of the Lie bialgebra} a definition of involutive Lie
bialgebra is given and in Appendix~\ref{curves}, we describe the Goldman-Turaev Lie
bialgebra of curves and we prove that it is involutive (\propc{involutive}).

In the final stages of the combinatorial treatment we benefited from the basic papers
\cite{BS} of Birman and Series, and \cite{CL} of Cohen and Lustig which helped us to
understand the relation between our combinatorics and hyperbolic geometric. Joint work
with Jane Gilman about a concrete homology intersection matrix for a surface with a
symmetric adapted basis \cite{G} also gave impulse and new ideas for   our efforts. This
work also benefited from discussions with Bill Goldman, Feng Luo,  Dennis Sullivan, and
Vladimir Turaev, and from a visit to Renaissance Technology.

\section{The vector space of cyclic words}\label{words}

\subsection{Cyclic words and linked pairs}

\begin{figure}[ht]
\begin{pspicture}(9,2)
\rput(7,1.2){$\curvearrowleft$}
 \rput(7,0.1){$a_1$}
\rput(7,1.9){$a_3$}\rput(7.9,1){$\ova_7$}\rput(6.1,1){$\ova_6$}\rput(7.67,1.57){$a_4$}
\rput(6.45,1.57){$a_5$} \rput(6.45,0.37){$\ova_4$} \rput(7.67,0.37){$a_4$}
\end{pspicture}
\caption{A ring in the letters of 
$\ensuremath{\mathbb{A}}_7$}\label{ring} 
\end{figure}
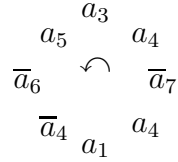

For each non-negative integer $n$, the {\em $n$-alphabet} or, 
briefly, the {\em alphabet} is the set of $2n$ symbols 
$\an=\conj{a_1,a_2\dots, 
a_n,\ov{a}_1,\ov{a_2},\dots,\ov{a}_n}.$  We shall consider 
linear words, denoted with capital roman characters, and 
cyclic words, denoted by capital caligraphycal characters, 
both in the letters of $\an$. The reader should think of 
cyclic words as symbols placed at the vertices of the $n$-th 
root of unity in $\mathbb{C}$   up to circular  symmetry, 
$n=1,2,3,\dots$, see Figure \ref{ring}.  

If $x_0x_1\dots x_{m-1}$ is a linear word, then, by 
definition, $\ov{x_0x_1\dots x_{m-1}}=\ox_{m-1}\ox_{m-2}\dots 
\ox_0$ and for each letter $x$, $\ov{\ox}=x$. A linear word 
$x_0x_1\dots x_{m-1}$ is {\em freely reduced} if $x_i\ne 
\ovx_{i+1}$ for each $i \in \conj{0,1,\dots,m-1}$. A linear 
word $\WWW_1$ is a {\em linear representative} of a cyclic 
word $\WW$ if $\WWW_1$ can be obtained from $\WW$ by making a 
cut between two consecutive letters of $\WW$. In such a case, 
we write $\WW=\cc(\WWW_1)$. If $\WW$ is a cyclic word, 
$\WWW_1$ is a linear representative of $\WW$, and $n$ is a 
positive integer,  we define $\WW^n$ as $\cc(\WWW_1^n)$, 
$\ov{\WW}$ as $\cc(\ov{\WWW}_1)$, and $\WW^{-n}$ as 
$\ov{\WW}^{\,n}$. (These are the basic operations on 
conjugacy classes in groups and are well defined by these 
prescriptions).  A cyclic word is {\em reduced} if it is 
non-empty and all its linear representatives are freely 
reduced, i.e., if the arrangement of symbols or ring is 
reduced. A reduced cyclic word is {\em primitive} if it 
cannot be written as $\WW^r$ for some $r \ge 2$ and some 
reduced cyclic word $\WW$. The {\em length} of a linear 
(resp. cyclic) word $\WWW$ (resp. $\WW$) is the number of 
letters counted with multiplicity that it contains and it is 
denoted by $\len(\WWW)$ (resp. $\len(\WW$)). By a {\em 
subword of a cyclic word $\WW$} we mean a linear subword of 
one of the linear representatives of $\WW$. 

Let $\OO$ be a  reduced cyclic word such that every letter of $\an$ appears exactly once.
From now on, we shall work with a fixed word $\OO$  and all our constructions will depend
on this choice (nevertheless, see \remc{isomorphic}). Such an $\OO$ is called the {\em
surface symbol}.

To each cyclic word $\WW$, we associate a number, 
$\si(\WW)\in \conj{-1,0,1}$ as follows. If $\WW$ is  reduced 
and there exists an injective orientation preserving (resp. 
orientation reversing) map of oriented rings, from the 
letters of $\WW$ to the letters of $\OO$ then $\si(\WW)=1$ 
(resp. $\si(\WW)=-1$). In all other cases (that is, if $\WW$ 
is not reduced or if there is no such orientation preserving 
or reversing map )   $\si(\WW)=0$.

\start{defi}{def} Let $\ap, \aq$ be two linear words. The ordered pair $(\ap,\aq)$ is
{\em $\OO$-linked} or briefly, {\em linked}, if $\ap$ and $\aq$ are reduced words of
length at least two and one of the following conditions holds.
\begin{numlist}
\item $\ap=p_1p_2$, $\aq=q_1q_2$ and  $\si(\cc(\ov{p}_1\ov{q}_1p_2q_2))\ne 0$.
\item $\ap=p_1\YYY  p_2$, $\aq=q_1\YYY q_2$, $p_1 \ne q_1$, $p_2 \ne q_2$ and $\YYY$ is a linear word of length at
least one, and if one writes $\YYY$ as $x_1\XXX x_2$ (where 
$\XXX$ is an empty word if $\YYY$ has length two and $x_1$ 
coincides with $x_2$ if the length of $\YYY$ is one), then  
$\si(\cc(\ov{p}_1 \ov{q}_1 x_1))=\si(\cc(p_2 q_2 \ox_2))$. 
\item $\ap=p_1\YYY  p_2$, $\aq=q_1\ov{\YYY} q_2$, $p_1 \ne \ov{q}_2$, $p_2 \ne \ov{q}_1$ and
$\YYY$ is a linear word of length at least one, and if one writes $\YYY$ as $x_1\XXX x_2$
(where $\XXX$ may be an empty word and $x_1$ may be equal to $x_2$), then
$\si(\cc(q_2\ov{p}_1x_1))=\si(\cc(\ov{q}_1p_2 \ox_2))$.
\end{numlist}
\end{defi}

Linked pairs capture the following idea: Two strands on a 
surface come close, stay together for some time and then 
separate. If one strand enters the strip from above and exits 
below and the other vice versa we must have an intersection. 
This is measured by linked pairs (see Figure \ref{motive}) 

The above definition surfaced studying the structure of the 
intersection and self-intersection points of curves in an 
orientable surface with boundary. It may seem obscure at the 
first reading but the reader will find the full motivation in 
Section~\ref{Intersection}. In particular, we will show how 
Figure~\ref{motive} illustrates this definition. 

If $\WW$ is a  reduced cyclic word, denote by $\bb_1 (\WW)$ 
the set of linked pairs $(\ap,\aq)$, where $\ap$ and $\aq$ 
are occurrences of  linear subwords of $\WW$. For 
understanding the definition better right now we suggest the 
reader do the calculations of Example~\ref{first}. 

\start{ex}{first} Set $\OO=\cc(a_1 a_2 \ova_1 \ova_2 a_3 a_4 
\ova_3 \ova_4)$  and consider $\WW=\cc(a_1 a_2 \ova_3 a_1 a_1 
a_3 \ova_2 a_1)$. There are linked pairs in $\bb_1 (\WW)$ of 
all types.
\begin{alphlist}
\item  the pairs; $\cc(a_3a_1a_3a_1a_2a_2)$ 
($1$) $(a_2\ova_3, \ova_3 a_1), (\ova_3 a_1,a_2\ova_3)$ 
$(a_2\ova_3, a_1 a_3),$  $(a_1 a_3, a_2\ova_3),$ $(\ova_3 
a_1, a_3 \ova_2),$ $(a_3 \ova_2,\ova_3 a_1),$ $(a_1 a_3, a_3 
\ova_2),$ $(a_3 \ova_2, a_1 a_3)$ satisfy 
\defc{def}(1);
\item  $(\ova_3 a_1 a_1, a_1 a_1 a_3),$ $(a_1 a_1 a_3, \ova_3 
a_1 a_1),$ $(\ova_2a_1 a_1, a_1 a_1 a_2),$ $(a_1 a_1 a_2, 
\ova_2a_1 a_1)$ satisfy \defc{def}(2); 
\item  
$(a_1a_2\ova_3a_1, a_1a_3\ova_2a_1),$ $(a_1a_3\ova_2a_1, 
a_1a_2\ova_3a_1)$ satisfy 
\defc{def}(3)
\end{alphlist}\end{ex}

\start{rem}{fin} Since a linear subword of a cyclic word is 
determined by its ordered pair of endpoints,  the set of 
occurrences of linear subwords in a  reduced cyclic word 
$\WW$ contains $\len(\WW)^2$ elements. Therefore, the set of 
linked pairs of $\WW$, $\bb_1(\WW)$ is finite and contains at 
most $\len(\WW)^4$ elements. A more careful study (\propc{is 
finite}) will show that $\len(\WW).(\len(\WW)-1)$ is an upper 
bound for the cardinality of $\bb_1(\WW)$. 
\end{rem}

\subsection{Definition of the cobracket}

We denote by $\V$ the vector space generated by non-empty 
reduced cyclic words in the letters of $\an$. Our next 
objective consists in defining a Lie cobracket on $\V$, which 
is a linear map $\map{\cob}{\V}{\V\otimes \V}$ satisfying 
certain identities (see Appendix \ref{definition of the Lie 
bialgebra}). Let us first motivate the definition: Observe 
that making two cuts between two different pairs of 
consecutive letters of a cyclic word, one gets two linear 
words. By gluing together the ends of each linear word, and 
reducing if necessary, one  obtains a pair of reduced  or 
empty cyclic words. On the other hand, as we shall soon see,  
every linked pair of subwords of a cyclic word $\WW$ 
determines two pairs of consecutive letters where one can 
make two cuts. Therefore, every linked pair of subwords of a 
cyclic word determines a pair of reduced or empty cyclic 
words. We will also see that these two cyclic words are 
non-empty (\propc{reduced cobracket}) and that the linked 
pair determines an ordering of this pair of cyclic words. 

Here is the precise definition of the procedure of the above paragraph. To each ordered
pair $(\ap,\aq) \in \bb_1 (\WW)$ we associate two cyclic words
$\cob_1(\ap,\aq)=\cc(\WWW_1)$ and $\cob_2(\ap,\aq)=\cc(\WWW_2)$ by the following.
\begin{romlist}
\item Assume that $(1)$ or $(2)$ of \defc{def} hold. Make two cuts on $\WW$,
 one immediately before $p_2$ and the other immediately before $q_2$. We obtain two linear words, $\WWW_1$
 and $\WWW_2$, the former, starting at $p_2$, and the latter, starting at $q_2$.
\item If condition $(3)$ holds, let $\WWW_1$ be the linear subword of $\WW$ starting at
$p_2$ and ending at $q_1$, and let $\WWW_2$ be the linear subword of $\WW$ starting at
$q_2$ and ending at $p_1$.
\end{romlist}

\start{lem}{disjoint} Let $\WW$ be a cyclic reduced word. For each $(\ap,\aq) \in
\bb_1(\WW)$, the linear words $\WWW_1$ and $\WWW_2$ of the above definition are disjoint
in $\WW$. Moreover, $\WWW_1$ and $\WWW_2$ are non-empty and one can write $\WW=\cc(\WWW_1
\WWW_2)$ in the case $(i)$ above and $\WW=\cc(\YYY\WWW_1\ov{\YYY}\WWW_2)$ in the case
$(ii)$, where $\YYY$ is as in \defc{def}(3).
\end{lem}
\begin{proof} The proof of the result in the case $(i)$ is straightforward. Let us study the
case $(ii)$. We claim that $\YYY$ and $\ov{\YYY}$ cannot overlap. Indeed, assume, for
instance that there exist a pair of linear words $\YYY_1$ and $\YYY_2$ such that
$\YYY=\YYY_1\YYY_2$ and $\ov{\YYY}$ starts with $\YYY_2$. Notice that, by definition,
$\ov{\YYY}=\ov{\YYY}_2\ov{\YYY}_1$. Therefore, $\YYY_2=\ov{\YYY}_2$. Since $\YYY_2$ is
reduced this cannot happen.

Finally, if one of the words, for instance, $\WWW_1$ is empty, then since $\YYY$ is
non-empty, $\WW=\cc(\YYY\ov{\YYY}\WWW_2)$ is not reduced, which contradicts our
hypothesis.

\end{proof}

\start{prop}{reduced cobracket} Let $\WW$ be a  reduced cyclic word and let $(\ap, \aq)
\in \bb_1(\WW)$ be a linked pair. Then $\cob_1(\ap,\aq)$ and $\cob_2(\ap,\aq)$ are
 reduced cyclic words. Moreover, $\cob_1(\ap,\aq)$ and $\cob_2(\ap,\aq)$
are non-empty.
\end{prop}
\begin{proof}Observe that, by
\lemc{disjoint}, the linear words  $\WWW_1$ and $\WWW_2$ are 
non-empty. Hence in order to prove that the cyclic words 
$\cob_1(\ap,\aq)$ and $\cob_2(\ap,\aq)$ are non-empty, it is 
enough to prove that   $\cc(\WWW_1)$ and $\cc(\WWW_2)$ are 
reduced. 

Assume that $(\ap,\aq)$ satisfies \defc{def}(1). Now, $\cob_1(\ap,\aq)$ is not reduced if
and only if $p_2=\ov{q}_1$ and  $\cob_2(\ap,\aq)$ is not reduced if and only if
$p_1=\ov{q}_2$. Since $\si(\cc(\ov{p}_1\ov{q}_1p_2q_2))\ne 0$, $p_2 \ne \ov{q}_1$ and
$p_1\ne \ov{q}_2$.

If $(\ap, \aq)$ satisfies \defc{def}(2) then $\cob_1(\ap,\aq)$ (resp. $\cob_2(\ap,\aq)$)
is not reduced if and only if $\ovx_2=p_2$ (resp. $\ovx_2=q_2$). Since $\WW$ is
 reduced, none of those equations can be satisfied.

Finally, if $(\ap, \aq)$ satisfies \defc{def}(3) then 
$\cob_1(\ap,\aq)$ is not reduced if and only if 
$p_2=\ov{q}_1$ and  $\cob_2(\ap,\aq)$ is not reduced if and 
only if $p_1=\ov{q}_2$. On the other hand, by \defc{def}(3), 
$p_2\ne \ov{q}_1$ and $p_1\ne \ov{q}_2$. Hence the proof is 
complete. 
\end{proof}

\start{ex}{more of first} Let $\OO$ and $\VV$ be as in 
\exc{first}. Hence for instance, $$\cob_1(a_2\ova_3, \ova_3 
a_1)=\cc(\ova_3) \mbox{ and } \cob_2(a_2\ova_3, \ova_3 
a_1)=\cc( a_1 a_1 a_3 \ova_2 a_1a_1 a_2),$$ 
$$\cob_1(a_1a_2\ova_3a_1, a_1a_3\ova_2a_1)=\cc(a_1a_1) \mbox{ 
and } \cob_2(a_1a_2\ova_3a_1, a_1a_3\ova_2a_1)=\cc(a_1 
a_1).$$ 
\end{ex}

By the following equation, one associates a sign to each linked pair  $(\ap, \aq)$.
\begin{eqnarray*}
\sign(\ap, \aq) &= & \left\{
        \begin{array}{ll}
\si(\cc(\ov{p}_1\ov{q}_1p_2q_2)) &\mbox{if $(\ap,\aq)$ satisfies  \defc{def}(1).}\\
\si(\cc(\ov{p}_1 \ov{q}_1 x_1)) &\mbox{if $(\ap,\aq)$ satisfies \defc{def}(2).}\\
\si(\cc(q_2\ov{p}_1x_1))&\mbox{if $(\ap,\aq)$ satisfies in \defc{def}(3).}
        \end{array}
\right.
\end{eqnarray*}

\start{lem}{minus} \begin{alphlist}
\item For every linked pair $(\ap, \aq)$, $\sign(\ap,\aq)=1$ or $\sign(\ap,\aq)=-1$.
\item If $(\ap, \aq)$ is a linked pair then
$(\aq,\ap)$ is also a linked pair. Moreover, $\sign(\ap,\aq)=-\sign(\aq,\ap)$.
\end{alphlist}
\end{lem}
\begin{proof} Let us prove $(a)$ and leave the proof of $(b)$ to the reader.
If $(\ap, \aq)$ satisfies \defc{def}(1), the result holds by definition. Otherwise,
observe that if $\ZZZ$ is a reduced cyclic word of three letters then $\si(\ZZZ)=1$ or
$\si(\ZZZ)=-1$, which implies the result.
\end{proof}

Now, using the above sign and the definitions of $\cob_1$ and $\cob_2$, we define
$\map{\cob}{\V}{\V\otimes \V}$ as the linear map such that for every  reduced cyclic word
$\WW$, $$\cob(\WW)=\sum_{(\ap,\aq)\in \bb_1(\WW)}\sign(\ap,\aq)\,\,\cob_1(\ap,\aq)\otimes
\cob_2(\ap,\aq).$$ By definition, the set  $\bb_1(\WW)$ is finite, hence, the above sum
is finite.

\subsection{Definition of the bracket}

Let $\VV$ and $\WW$ be two cyclic words and choose a pair of 
consecutive letters in each word. Performing a cut between 
each of these pair of letters, one obtains two linear words, 
$\VVV_1$ and $\WWW_1$. The linear word $\VVV_1 \WWW_1$ 
determines a cyclic word, $\cc(\VVV_1 \WWW_1)$ (possibly not  
reduced). In other words, a pair of cuts on a pair of cyclic 
words determines a third cyclic word. Now, we want to define 
an object analogous to $\bb_1(\WW)$ not for a single cyclic 
word $\WW$ but for a pair of cyclic words $\VV$ and $\WW$.  
This object will be denoted by $\bb_2(\VV,\WW)$ and should 
contain linked pairs $(\ap, \aq)$ such that $\ap$ is a 
subword of $\VV$ and $\aq$ is a subword of $\WW$ (and, as we 
shall see, it must also contain other linked pairs). Since 
the formal definition could seem unnatural, we have chosen to 
give first a motivation for it. The first naive definition is 
(at least, so it was for the author of this paper) let the 
set of linked pairs of a pair of cyclic words $\VVV$ and 
$\WWW$ be the set of linked pairs $(\ap, \aq)$ such that 
$\ap$ is a subword of $\VVV$ and $\aq$ is a subword of 
$\WWW$. The following example illustrates the limitations of 
this naive definition.

\start{ex}{problem} Let us consider the alphabet 
$\atwo=\{a_1,a_2, \ova_1,\ova_2\}$ and define $\OO=a_1 a_2 
\ova_1 \ova_2$. Set $\VV=a_1 a_1 a_2$, $\WW=a_1 a_1 a_2a_1 
a_1 a_2a_1$, $\ap=a_2 a_1 a_1 a_2a_1 a_1 a_2a_1 a_1 a_2$ and 
$\aq=a_1 a_1 a_1 a_2a_1 a_1 a_2a_1 a_1 a_1$. Then $(\ap, 
\aq)$ is an $\OO$-linked pair. On the other hand, $\ap$ is  a 
subword of $\VV^4$ but not a subword of $\VV$, $\VV^2$ or 
$\VV^3$. Analogously, $\aq$ is a subword of $\WW^2$ but not a 
subword of $\WW$. 
\end{ex}

Because linked pairs correspond to intersection points we 
will have to  consider subwords of all powers of pairs of 
cyclic words as the set of linked pairs of a two cyclic 
words, (see below). Unlike the case of the set of linked 
pairs of a single word $\bb_1(\WW)$, the definition of 
$\bb_2(\VV,\WW)$ does not obviously imply that this set is 
finite. The finiteness of $\bb_2(\VV,\WW)$ is true and the 
size will be estimated in \propc{is finite}. Another issue 
concerning linked pairs of subwords of all powers of a cyclic 
word is whether the length of the subwords is bounded. It 
turns out that there exists a bound on the length, depending 
on the length of the pair of original cyclic words 
(\propc{finite length}). Notice that \propc{finite length} 
also implies that the cardinality of the set of linked pairs 
of two words, $\bb_2(\VV,\WW)$ is bounded, but unlike 
\propc{is finite} it does not give a sharp upper bound for 
the number of elements in $\bb_2(\VV,\WW)$.

For each pair of cyclic reduced words, $\VV$ and $\WW$, the 
{\em set of linked pairs of $\VV$ and $\WW$}, denoted by 
$\bb_2(\VV,\WW)$, is defined to be the set of all pairs 
$(\ap, \aq)$ for which there exists positive integers $j$ and 
$k$ such that $\ap$ is an occurrence of a subword of $\VV^j$, 
$\aq$ is an occurrence of a subword of $\WW^k$, where
$\len(\VV^{j-1})<\len(\ap)\le \len(\VV^{j})$ and 
$\len(\WW^{k-1})<\len(\aq)\le \len(\WW^{k})$. (Here, we set 
$\len(\VV^0)=\len(\WW^0)=0$.)

\start{prop}{is finite} Let $\VV$ and $\WW$ be reduced cyclic words. Then
\begin{alphlist} \item There are at most $\len(\VV).\len(\WW)$ elements in
$\bb_2(\VV,\WW)$, the set of linked pairs of  $\VV$ and $\WW$.
\item The set of linked pairs of one word, $\bb_1(\WW)$ contains at most
$\len(\WW).(\len(\WW)-1)$ elements.
\end{alphlist}
\end{prop}
\begin{proof} We first prove $(a)$. Let $n_{1}, n_2, n_3$ be the number of linked pairs
satisfying $(1)$, $(2)$ and $(3)$ respectively of \defc{def} for a pair of cyclic reduced
words $\VV$ and $\WW$.

Let $C$ be the set of pairs of the form $(x_1x_2, z_1 z_2)$ such that $x_1x_2$ is an
occurrence of a subword of $\VV$, and $z_1z_2$ is an occurrence of a subword of $\WW$.
Notice that the cardinality of $C$ is $\len(\VV).\len(\WW)$.

Let $C_1$ be  the set of pairs $(x_1x_2, z_1 z_2) \in C$ such that $x_1\ne z_1$ and $x_2
\ne z_2$. Clearly, the set of linked pairs satisfying \defc{def}$(1)$ is contained in
$C_1$ and so any upper bound of the cardinality of $C_1$ is larger than $n_1$.

We claim that each linked pair satisfying \defc{def}$(2)$ 
determines two different elements in $C \setminus C_1$, and 
that for every positive integer $k$, $k$ different linked 
pairs satisfying \defc{def}$(2)$ determine $2k$ such 
elements. By this claim, since there are $n_2$ pairs 
satisfying \defc{def}$(2)$, there are $2n_2$ pairs in 
$C\setminus C_1$. Then, the cardinality of $C_1$ is at most 
$\len(\VV).\len(\WW)-2n_2$ and so 
$n_1\le\len(\VV).\len(\WW)-2n_2$. 

In order to prove the claim, consider a linked pair satisfying $(2)$, let $\YYY$ be the
``middle'' linear word of this pair as in \defc{def}(2) and let $y$ be the first letter
of $\YYY$. Let $(\ap_1,\aq_1)$, $(\ap_2,\aq_2)$ by the elements of $C$, such that $\ap_1$
and $\aq_1$ start with the first letter of $\YYY$ and $\ap_2$ and $\aq_2$ end with the
first letter of $\YYY$. Since all the linear words $\YYY$ are different and, by
definition,  $(\ap_1,\aq_1)$ and $(\ap_2,\aq_2)$ are not in $C_1$, the claim is proved

Let $C_2$ be  the set of pairs $(x_1x_2, z_1 z_2) \in C$ such 
that $\ovx_1\ne z_2$ and $\ovx_2 \ne z_1$. The set of linked 
pairs satisfying $(1)$ is contained in $C_2$. By an argument 
similar to the one we used to determine the cardinality of 
$C_1$, one can show that the cardinality of $C_2$ is 
$\len(\VV).\len(\WW)-2n_3$. Hence 
$n_1\le\len(\VV).\len(\WW)-2n_3$. 

Adding both inequalities, one gets $2n_1\le 
2\len(\VV).\len(\WW)-2n_2-2n_3$. Hence $n_1+n_2+n_3\le 
\len(\VV).\len(\WW)$, and the proof of $(a)$ is complete. 

In order  to prove $(b)$ we proceed analogously as we did in $(a)$, but defining $C$ as
the set of pairs of linear two-letter subwords $(\ap, \aq)$ of $\VV$ such that $\ap$ and
$\aq$ do not start at the same letter, i.e., $\ap$ and $\aq$ are different occurrences of
subwords of $\VV$. Hence $C$ has $\len(\WW).(\len(\WW)-1)$ elements and all the pairs
satisfying \defc{def}(1) are in $C$. Now, we can complete the proof as above.

\end{proof}

\start{ex}{sharp} The following examples show that the bounds 
of \propc{is finite} are sharp. Let 
$\OO=a_1a_2\ova_1\ova_2a_3a_4\ova_3\ova_4$. Then 
$\bb_1(\cc(a_1\ova_3a_2))$ contains exactly six elements, 
which are $$\{ (a_1\ova_3,\ova_3a_2), (a_1\ova_3,a_2a_1), 
(\ova_3a_2,a_1\ova_3), (\ova_3a_2,a_2a_1), 
(a_2a_1,a_1\ova_3),(a_2a_1,\ova_3a_2) \}.$$ Also 
$\bb_2(\cc(a_1\ova_3),\cc(a_2\ova_4))$ contains four elements 
which are, $$\conj{(a_1\ova_3, a_2\ova_4),(a_1\ova_3, 
\ova_4a_2),(\ova_3a_1, a_2\ova_4),(\ova_3a_1, \ova_4a_2)} $$ 
\end{ex}

The next result about linear words is well known and it will be used in the proofs of the
Lemmas \ref{subwords} and \ref{pair}.

\start{lem}{well} If $\ap=x_0x_1\dots x_{m-1}$ is a linear word and for some $i \in
\conj{1,2,\dots,m-1}$,
 $\ap=x_{i}x_{i+1}\dots x_{m-1}x_0 x_1\dots x_{i-1}$ then there exists a linear word $\aq$ and an integer $r$ such that
$r \ge 2$, $\ap=\aq^r$ and $\len(\aq)$ divides $i$.
\end{lem}

The next lemma will be used in the proof of \propc{finite length}.

\start{lem}{subwords} Let $\VV$, $\WW$ be cyclic words which 
are not powers of the same cyclic word and let $\ap$ be a 
linear word. Let $k, l$ be a pair of positive integers such 
that $\ap$  is  a subword of $ \VV^k$ and either $\ap$ or 
$\ov{\ap}$ is a subword of $\WW^l$. Moreover, assume that 
$(k-1)\len(\VV)<\len(\ap)$ and $(l-1)\len(\WW)<\len(\ap)$. 
Then  $\len(\ap)<\len(\VV)+\len(\WW)$. 

\end{lem}
\begin{proof} Let us  first consider the case that $\ap$ is a subword of both, $\VV^k$ and $\WW^l$.
Assume that $\len(\ap)\ge \len(\VV)+\len(\WW)$. Then we can 
write $\ap=p_0p_1p_2\dots p_{m-1}$ where $m 
\ge\len(\VV)+\len(\WW)$. Since $\ap$ is a subword of $\WW^l$, 
$$p_0p_1\dots p_{\len(\VV)-1}=p_{\len(\WW)}p_{\len(\WW)+1} 
\dots p_{\len(\WW)+\len(\VV)-1}.$$ Since $\ap$ is a subword 
of $\VV^k$, $$p_{\len(\WW)}p_{\len(\WW)+1}\dots 
p_{\len(\WW)+\len(\VV)-1}=p_rp_{r+1}\dots 
p_{r+\len(\VV)-1}=p_rp_{r+1}\dots p_{\len(\VV)-1}p_0p_1\dots 
p_{r-1},$$ where $r$ is the remainder of dividing $\len(\WW)$ 
by $\len(\VV)$. Since $\VV$ is not a power of $\WW$, $r>0$. 
By \lemc{well}, there exists a linear word $\XXX$ and a 
positive integer $d$ such that $\len(\XXX)$ divides $r$ and 
$p_0p_1\dots p_{\len(\VV)-1}=\XXX^d$. So $\VV=\cc(\XXX)^d$. 
Since $\len(\XXX)$ divides $r$ and $\len(\VV)$, $\len(\XXX)$ 
divides $\len(\WW)$. Thus, $\WW$ is also a power of 
$\cc(\XXX)$, contradicting our hypothesis. Therefore, 
$\len(\ap)<\len(\VV)+\len(\WW)$.

To prove the result in the other case, one observes that $\ov{\ap}$ is a subword of a
cyclic word $\WW$ if and only if $\ap$ is a subword of $\ov{\WW}$.  Then replacing  $\WW$
by $\ov{\WW}$, the result follows.
\end{proof}

Up to certain special cases left to the reader, the following proposition follows from
\lemc{subwords}.

\start{prop}{finite length} Let $\VV$ and $\WW$ be two 
reduced cyclic words. Then $\bb_2(\VV,\WW)$ is the set of all 
linked pairs $(\ap, \aq)$ such that $\ap$ is an occurrence of 
a subword of $\VV^j$, $\aq$ is an occurrence of a subword of 
$\WW^k$, $\len(\VV^{j-1})<\len(\ap)\le \len(\VV^{j})$ so that  
$\len(\WW^{k-1})<\len(\aq)\le \len(\WW^{k})$ where $j,k$ are 
positive integers such that $j<2+\frac{\len(\WW)}{\len(\VV)}$ 
and $k<2+\frac{\len(\VV)}{\len(\WW)}$. 
\end{prop}

Now, we will define the bracket. Firstly, we associate to 
each linked pair $(\ap,\aq)\in \bb_2(\WW,\ZZ)$ a cyclic word 
$\br(\ap,\aq)=\cc(\WWW_1\ZZZ_1)$, where $\WWW_1$ and $\ZZZ_1$ 
are linear words defined as follows. 
\begin{romlist}
\item If conditions $(1)$ or $(2)$ of \defc{def} holds for the linked pair
$(\ap,\aq)$, $\WWW_1$ is the representative of $\WW$ obtained by cutting $\WW$
immediately before $p_2$ and $\ZZZ_1$ is the representative of $\ZZ$ obtained by cutting
$\ZZ$ immediately before $q_2$.
\item If condition $(3)$ of \defc{def} holds for the pair $(\ap,\aq)$. Then $\WWW_1$ is the
linear subword of $\WW$ that starts right after the end of $\YYY$ and ends right before
the first letter of $\YYY$, and $\ZZZ_1$ is the subword of $\ZZ$ that starts right after
the last letter  of $\ov{\YYY}$ and ends right before the beginning of $\ov{\YYY}$.
(Observe that $\YYY$ may not be a subword of $\WW$ nor of $\ZZ$, but we can always find
the first and last letters of $\YYY$  in $\WW$ and $\ZZ$)
\end{romlist}

\start{ex}{c} Set $\OO=\cc(a_1 a_2 \ova_1 \ova_2 a_3 a_4 \ova_3 \ova_4)$. Then
$$\bb_2(\cc(a_1 a_2 a_2 a_3),\cc(\ova_2 \ova_2))=\conj{(a_1 a_2 a_2 a_3,\ova_2
\ova_2\ova_2 \ova_2), (a_1 a_2 a_2 a_3,\ova_2 \ova_2\ova_2 \ova_2)},$$ and both pairs
$(a_1 a_2 a_2 a_3,\ova_2 \ova_2\ova_2 \ova_2)$, $(a_1 a_2 a_2 a_3,\ova_2 \ova_2\ova_2
\ova_2)$ satisfy $(3)$ of \defc{def}. Notice that $\ZZZ_1$ in this case is the empty
word.
\end{ex}

The next result is  analogous to \propc{reduced cobracket}.

\start{prop}{reduced bracket} For each pair of reduced cyclic words $\WW$ and $\ZZ$, and
for each linked pair  $(\ap, \aq) \in \bb_2(\WW,\ZZ)$, $\br(\ap,\aq)$ is a cyclically
reduced word. In particular, $\br(\ap,\aq)$ is non-empty.
\end{prop}
\begin{proof} This proof follows the same ideas of the proof of \propc{reduced
cobracket}, except that in the case $(ii)$, the words 
$\WWW_1$ or $\ZZZ_1$ may be empty. Hence it is necessary to 
see that they cannot both be empty. Observe that if $\WWW_1$ 
and $\ZZZ_1$ are  empty, then $\YYY=\WWW_2^n$, 
$\ov{\YYY}=\ZZZ_2^m$, for some positive integers $n, m$, and  
some linear representative $\WWW_2$ of $\WW$ and $\ZZZ_2$ of 
$\ZZ$. But in this case, if we write $\ap=p_1\YYY p_2$, 
$\aq=q_1\ov{\YYY} q_2$,  as in 
\defc{def}, we have that $p_1 = \ov{q}_2$, $p_2 =\ov{q}_1$. So $(\ap, \aq)$ does not
satisfy $(3)$ of \defc{def}.\end{proof}

We define the bracket $\map{\bra{,}}{\V\otimes \V}{\V}$ as the linear map such that for
each pair of cyclic words, $\WW$ and $\ZZ$, $$\bra{\WW,\ZZ}=\sum_{(\ap,\aq)\in
\bb_2(\WW,\ZZ)}\sign(\ap,\aq)\,\,\br(\ap,\aq).$$

By \propc{is finite}, $\bb_2(\VV,\WW)$ is finite and so the bracket is well defined.

\start{ex}{d} If $\OO=\cc(a_1 a_2 \ova_1 \ova_2 a_3 a_4 \ova_3 \ova_4)$  then $[\cc(a_1
a_2 a_2 a_3),\cc(\ova_2 \ova_2)]=-2 \cc(a_3 a_1).$ (see \exc{c})
\end{ex}

In Section \ref{GT} we prove that for each $\OO$, $(\V, \br, \cob)$ is an involutive Lie
bialgebra (\theoc{Main}).

\start{rem}{isomorphic} Even though there is a vast number of 
$\OO$'s for each $n$ (exactly, $(n-1)!$), and  each $\OO$ 
determines a  Lie bialgebra,  there are at most $n/2+1$ Lie 
bialgebras up to isomorphism, (since a surface is determined 
by its genus and the number of boundary components.)\end{rem}

\section{Intersection points and linked pairs}\label{Intersection}

The goal of this section is to prove that there are bijective  
correspondences between ``unordered'' linked pairs of a 
cyclic reduced word and self-intersection points of 
representatives with minimal self-intersection and between 
linked pairs of a pair of words and intersection points of a 
pair of  representatives  with minimal intersection. Our 
representatives are similar to the ones constructed by 
Reinhart \cite{R} but improves them in the sense that 
intersection and self-intersections of our representatives 
are minimal and are organized by the combinatorics of linked 
pairs. 

We start by gathering  together some well known results.

\start{lem}{correspondence} Let $\Sigma$ be an oriented surface with boundary. There is a
bijective correspondence between any two of the following sets.
\begin{numlist}
\item   Conjugacy classes of non-trivial elements of the fundamental group of $\Sigma$.
\item  Non-trivial free homotopy classes of maps from the circle to $\Sigma$.
\item Non-empty cyclically reduced cyclic words in $\an=\conj{a_1,a_2\dots,
a_n,\ov{a}_1,\ov{a_2},\dots,\ov{a}_n},$ where $1-n$ is the Euler characteristic of
$\Sigma$.
\end{numlist}
\end{lem}

\lemc{correspondence} allows us to identify non-empty 
cyclically reduced cyclic words and non-trivial free homotopy 
classes, and we will often make use of this identification.

\subsection{Arc representatives}

As in Section~\ref{words}, throughout this section  we fix a surface symbol, that is, a
cyclic word $\OO=\cc(o_1o_2\dots o_{2n})$ such that every letter of $\an$ appears exactly
once. Denote by $\Po$ the $4n$-gon with edges labeled counterclockwise in the following
way: One chooses an edge as first and labels it with $o_1$, the second has no label, the
third is labeled  with $o_2$, the fourth with no label and so on as is shown in the
example of \figc{$4n$-gon}.

\begin{figure}
\begin{pspicture}(14,4)
\rput*{0}(12,2){$\OO=\cc(a_1 \ova_1 a_2 a_3 \ova_2 \ova_3)$}
\rput*(7,2){$\circlearrowleft$} \rput{14}(7,2){\PstPolygon[unit=1.6,PolyNbSides=12]}
\rput{0}(8.9,1.95){$a_1$} \rput(8,3.55){$\ova_1$}\rput(6,3.55){$a_2$}
\rput{0}(5.2,2){$a_3$}\rput{0}(6.05,0.55){$\ova_2$}\rput(8.1,0.55){$\ova_3$}
\end{pspicture}
\caption{The $4n$-gon $\Po$}\label{$4n$-gon}
\end{figure}

For each $i \in \conj{1,2,\dots,n}$, one identifies the edge $a_i$ with the edge
$\ov{a}_i$ without creating Moebius bands. In this way, one gets a surface $\So$ with
non-empty boundary and Euler characteristic $1-n$. Furthermore, every surface with
non-empty boundary can be obtained from such a $4n$-gon (one has to take $\OO$ as the
empty word in order to get the disk). Denote by $\map{\pi}{\Po}{\So}$ the projection map
and denote by $E_{x_i}$ the projection of the edge labeled with $x_i$. Observe that
$E_{x_i}=E_{\ovx_i}$.

A {\em loop in $\So$} is a piecewise smooth map from the circle to $\So$. We will often
identify a loop with its image.

Choose a point $c$ in the interior of $\Po$. For each $i \in \conj{1,2,\dots, n}$, let
 $q_i$ and  $s_i$ be a pair of points such that  
$\pi(q_i)=\pi(s_i)$, and that $q_i$ (resp. $s_i$) is in the 
interior of the edge labeled with $a_i$ (resp. $\ova_i$). 
Denote by $a_i$ the based homotopy class of the loop that 
starts at $\pi(c)$, runs along the projection of the segment 
from $c$ to $q_i$ and then along the projection of the 
segment  from $s_i$ to $c$.  Hence $\conj{a_1, a_2, \dots, 
a_n}$ is a set of generators of the fundamental group of 
$\So$. If $\WW$ is a reduced cyclic word in the letters of 
$\an$, a loop $\alpha$ in $\So$ is a {\em representative of 
$\WW$} if $\alpha$ is freely homotopic to a curve $\beta$ 
passing through $c$ such that the (based) homotopy class of 
$\beta$ written in the generators $\conj{a_1, a_2, \dots, 
a_n}$  is a linear representative of $\WW$ (see 
\lemc{correspondence}).

A one dimensional manifold (possibly with boundary)  $\beta$ (resp. a pair of loops
$\alpha, \beta$) on a surface $\Sigma$ {\em has a bigon} if there exists a pair $u$ and
$v$ of subarcs of $\beta$ (resp. such that $u$ is a subarc of $\alpha$ and $v$ is a
subarc of $\beta$), such that endpoints of $u$ equal the endpoints of $v$ and the loop
formed by running first along $u$ and then along $v$ is null-homotopic in $\Sigma$.

From now on, when considering a cyclic word $\WW=\cc(x_0x_1x_2\dots x_{m-1})$, for each
$i \in \Z$,  $x_i$ denotes the letter of $\WW$ with subindex $i$ modulo $m$. This
convention also applies to other objects with subindex.

A \emph{segment} is a map from a closed interval into the 
surface. If $A_0, A_1, A_2,\dots, A_{s-1}$ is a sequence of 
oriented segments  of $\Po$ such that for each $i \in 
\conj{0,1,\dots, s-2}$, the projection of the final point of 
$A_i$ is equal to the projection of the initial point of 
$A_{i+1}$, then $\pi(A_0)\pi(A_1)\dots \pi(A_{s-1})$ is an 
arc in $\So$. We will write $\pi(A_i A_{i+1}\dots A_{i+j})$ 
instead of $\pi(A_i)\pi(A_{i+1})\dots \pi(A_{i+j})$. Clearly, 
if the endpoint of $\pi(A_{s-1})$ equals the initial point of 
$\pi(A_0)$ then $\pi(A_0 A_{1}\dots A_{s-1})$ is a loop.

 Let $\WW=\cc(x_0x_1x_2\dots x_{m-1})$ be a
reduced cyclic word and let $\alpha$ be a representative of $\WW$. We say that $\alpha$
is a {\em segmented representative of $\WW$} if there exist a sequence of oriented
segments $A_0, A_1, A_2,\dots, A_{m-1}$ in $\Po$ such that $\alpha=\pi(A_0A_1\dots
A_{m-1})$, and
\begin{numlist}
\item For each $i \in \conj{0,1,2,\dots,m-1}$, $A_i$
is an oriented arc starting at a point in the interior of the edge labeled with $\ovx_i$,
and ending at a point  in  the interior of the edge labeled with $x_{i+1}$.
\item For each $i,j \in \conj{0,1,2,\dots,m-1}$, $A_i$ intersects $A_j$ in at most one
point.
\item The endpoints of the arcs $A_0$, $A_1,\dots, A_{m-1}$, are all different.
\item There are no triple intersections between the arcs $A_0, A_1, \dots, A_{m-1}$.
\item For each $i \in \conj{0,1,2,\dots,m-1}$, the  interior of $A_i$ is contained in the interior of $\Po$
\end{numlist}

The following lemma   characterizes  the segmented 
representatives of a reduced cyclic word with bigons. 

\start{lem}{arc bigon} Let $\WW=\cc(x_0x_1x_2\dots x_{m-1})$ be a  reduced cyclic word
and let $\alpha=\pi(A_0A_1\dots A_{m-1})$ be a segmented representative of $\WW$. Then
$\alpha$ has a bigon if and only if there exist $i, j \in \conj{0,1,\dots, m-1}$, $j \ne
i$, and $k \ge 0$ such that one of the following holds
\begin{numlist}
\item $\pi(A_i) \cap \pi(A_j) \ne \emptyset$, $\pi(A_{i+k}) \cap \pi(A_{j+k}) \ne
\emptyset$ and for each $h \in \conj{1,2,\dots, k-1}$,
$\pi(A_{i+h})\cap\pi(A_{j+h})=\emptyset$.
\item $\pi(A_{i}) \cap \pi(A_{j}) \ne \emptyset$, $\pi(A_{i+k}) \cap \pi(A_{j-k}) \ne
\emptyset$ and for each $h \in \conj{1,2,\dots, k-1}$,
$\pi(A_{i+h})\cap\pi(A_{j-h})=\emptyset$.
\end{numlist}
\end{lem}
\begin{proof} Clearly, if $(1)$ or $(2)$ hold then $\alpha$ has a bigon.  We prove now
the reverse implication. Let $U$, $V$ be the subarcs that bound a bigon, and let $p, q$
be the endpoints of $U$ and $V$. Assume that $U$ goes from $p$ to $q$ and $V$ goes from
$q$ to $p$. The proof in the other case is similar.

Let $i, j \in \conj{0,1,\dots, m-1}$, $s, l\ge 0$, be such 
that $p \in \pi(A_i) \cap \pi(A_j)$,  $q \in \pi(A_{i+s}) 
\cap \pi(A_{j-l})$, and $s$ and  $l$ are minimal with such a 
property. Let $\ap=x_{i+1}x_{i+2}\dots x_{i+s}$ and 
$\aq=x_{j-l}x_{j-l+1}\dots x_{j-1}$. Since the union of $U$ 
and $V$ is a null-homotopic loop, and $\WW$ is reduced, 
$\ap=\ov{\aq}$. Hence $s=l$ and taking $h=s$, the integers 
$i, j$ and $h$ satisfy $(2)$. 
\end{proof}

\start{prop}{gon} Let $\WW=\cc(x_0x_1x_2\dots x_{m-1})$ be a primitive reduced cyclic
word. Then there exist a segmented representative $\alpha=\pi(A_0A_1\dots A_{m-1})$ of
$\WW$ such that $\alpha$ does not have bigons.
\end{prop}

\begin{proof} Consider $U_i$ and $V_i$ thin tubular neighborhoods of the edges $a_i$ and $\ova_i$ in $\Po$
respectively. The projection of $U_i$ and $V_i$ determine two sides of $E_i$: the
$a_i$-side, containing $\pi(U_i)$  and the $\ova_i$-side, containing $\pi(V_i)$.

Now, to each loop $\varphi$ such that each small arc of 
$\varphi$ intersects  $\cup_{1\le i \le n}E_{a_i}$ 
transversely, we associate a cyclic word $\VV_\varphi$, which 
describes the way that $\varphi$ crosses $\cup_{1\le i \le 
n}E_{a_i}$: Choose a point $p$ in $\varphi$ and not in 
$\cup_{1\le i \le n}E_{a_i}$. Let $y_0y_1\dots y_{k-1}$ be an 
ordered sequence of letters of $\an$ such that, starting at 
$p$, the first edge  crossed by $\varphi$ is $E_{y_0}$, from 
the $y_0$-side to the $\ovy_0$-side. The second edge crossed 
by $\varphi$ is $E_{y_1}$ from the $y_1$-side to the 
$\ovy_1$-side and, in general, for each $i \in 
\conj{1,\dots,k}$, the $i$-th edge crossed by $\varphi$  is 
$E_{y_{i-1}}$ from the $y_{i-1}$-th side to the 
$\ovy_{i-1}$-side. Finally, take $\VV_\varphi=\cc(y_0y_1\dots 
y_{k-1})$.

Endow $\So$ with a hyperbolic metric with geodesic boundary  
such that the arcs $\cup_{1\le i \le n}E_{a_i}$ are also 
geodesic. (Such a metric exists because we can assume that 
$\Po$ is a hyperbolic polygon $\Po$ with geodesic edges and 
right angles.) Let $\beta$ be a geodesic representative of 
$\WW$. Then, $\beta \cup(\cup_{1\le i \le n}E_{a_i})$ does 
not have bigons, because a bigon cannot be bounded by 
geodesic segments. 

Let us prove that $\VV_\beta$ is reduced and therefore, 
$\VV_\beta=\cc(x_0x_1\dots x_{m-1})$. Indeed, suppose that 
$\VV_\beta=\cc(y_0y_1\dots y_{k-1})$ is not reduced. Hence 
$y_{i+1}=\ovy_i$, for some $i \in \conj{0,1,2,\dots, k-1}$. 
Then there exists a subarc $B$ of $\beta$, a subarc $S$ of 
$E_{y_i}$ such that $S$ and $B$  bound a disk, which is 
absurd, being a bigon bounded by geodesics.

Now, we claim that there exists an homotopy $\beta_t$, $t \in [0,1]$ such that:
\begin{romlist}
\item $\beta_0=\beta$.
\item For each $t \in [0,1]$, $\beta_t \cup(\cup_{1\le i \le n}E_{a_i})$
 does not have bigons.
\item  All the self-intersection points of $\beta_1 \cup(\cup_{1\le i \le n}E_{a_i})$ are
double.
\item All the  intersection points of pairs of arcs of
$\beta_t \cup(\cup_{1\le i \le n}E_{a_i})$ are transverse.
\end{romlist}

\begin{figure}[ht]
\begin{pspicture}(15.5,3)
 \rput(1,0){\rput(0.6,2.5){$\DD$}\pscircle(1.5,1.5){1}}
\psline[linecolor=darkgray](1.8,2.8)(2,2.4)(3,0.6)(3.2,0.3) \rput(2.9,0.2){$E_{a_i}$}
\pscurve(1.2,1.4)(1.5,1.5)(2.5,1.5)(3.4,1)(3.6,0.8)\pscurve(1.2,0.8)(1.6,1)(2.5,1.5)(3,2.4)(3.2,2.8)
\pscurve(1.8,0.4)(2,0.7)(2.5,1.5)(2.4,2.5)(2.2,2.8)
\rput(4,0){\rput(2.6,2.5){$\DD$} \pscircle(1.5,1.5){1}}
\rput(3,0){\psline[linecolor=darkgray](1.8,2.8)(2,2.4)(3,0.6)(3.2,0.3)
\rput(2.9,0.2){$E_{a_i}$}

\pscurve(1.2,1.4)(1.5,1.5)(2.5,1.5)(3.4,1)(3.6,0.8)
\pscurve(1.2,0.8)(1.6,1)(2.5,1.5)(3,2.4)(3.2,2.8)
\pscurve(1.8,0.4)(2,0.7)(2.5,1.5)(2.4,2.5)(2.2,2.8)

\pscircle[fillcolor=white, fillstyle=solid](2.5,1.5){1}
 \pscurve(1.5,1.5)(3,1.6)(3.4,1)
\pscurve(1.6,1)(3,1.6)(3,2.4) \pscurve(2,0.7)(3,1.6)(2.4,2.5)
\psline[linecolor=darkgray](1.8,2.8)(2,2.4)(3,0.6)(3.2,0.3)}

\rput(4,2.5){$\beta_t$} \pscurve{->}(3.6,2)(4,2.2)(4.4,2)
\rput(7,0){\rput(4,2.5){$\beta_t$}

 \rput(6,0){\rput(0.6,2.5){$\DD$} }\pscurve{->}(3.6,2)(4,2.2)(4.4,2)}

\pscircle(9.5,1.5){1} \rput(7,0){ \pscurve(1.2,1.4)(1.5,1.5)(2.5,1.5)(3.4,1)(3.6,0.8)
\pscurve(1.2,0.8)(1.6,1)(2.5,1.5)(3,2.4)(3.2,2.8)
\pscurve(1.8,0.4)(2,0.7)(2.5,1.5)(2.4,2.5)(2.2,2.8) }

\rput(10,0){ \pscurve(1.2,1.4)(1.5,1.5)(2.5,1.1)(3.4,1)(3.6,0.8)
\pscurve(1.2,0.8)(1.6,1)(2.5,2)(3,2.4)(3.2,2.8)
\pscurve(1.8,0.4)(2,0.7)(2.4,1.5)(2.4,2.5)(2.2,2.8) } \rput(6,0){\rput(2.6,2.5){$\DD$}}

\rput(4,0.2){$(a)$}\rput(11,0.2){$(b)$}

\pscircle(12.5,1.5){1}
\psdots[dotstyle=o,dotsize=2.2pt](8.6,1.2)(8.8,1.2)(9,1.2)
(8.6,1.4)(8.8,1.4)(9,1.4)(9.4,0.8)(9.6,0.8)(9.8,0.8)(10,0.8)(10,1)
(9.4,1)(9.6,1)(9.8,1)(10.2,1) (9.6,1.4)
(11.6,1.4)(11.8,1.4)(11.6,1.23)(11.8,1.23)(13.16,0.94)(12.6,2.2)(12.6,2.4)(12.8,2.4)
(12.76,2.22)
(9.2,1.4)(9.6,1.8)(9.6,2)(9.6,2.2)(9.6,2.4)(9.8,2.3)(9.77,2.1)(9.2,0.6)(9.2,0.8)(9.4,0.6)
(9.8,0.6)(9.6,0.6)(9.6,0.8)(9.6,1)(9.4,1.2)(9.6,1.2)(9.8,1.2) (10,1.2)
(12.2,0.6)(12.4,0.6)(12.6,0.6)(12.8,0.6)
(12.2,0.8)(12.4,0.8)(12.6,0.8)(12.8,0.8)(13,0.8)(12.4,1)(12.6,1)(12.8,1)(13,1)

\end{pspicture} \caption{Proof of \propc{gon}}\label{beta}
\end{figure}

Indeed, if  all the self-intersection points of $\beta\cup 
(\cup_{1\le i \le n}E_{a_i})$ are double, the claim holds 
trivially. Otherwise, we remove  self-intersection points 
which are not double in the following way: Let $p$ be one 
such point.  Consider a small disk $\DD$ centered  at $p$ 
such that  each connected component of $\DD\cap (\beta 
\cup(\cup_{1\le i \le n}E_{a_i}))$ passes through $p$. First, 
assume that one of these connected components is included in  
$\cup_{1\le i \le n}E_{a_i}$. In this case, deform $\beta$ as 
in Figure \ref{beta}(a) and then as in Figure \ref{beta}(b). 
Observe that in the second step, if $h$ is the number of 
connected components of $\beta \cap \DD$, then after the 
homotopy there are exactly $\frac{h(h-1)}{2}$ transverse 
double points in $\DD$. We can assume that in both 
homotopies, $\beta \setminus \DD$ is fixed. Now, assume that 
none of the connected components of $\DD\cap (\beta 
\cup(\cup_{1\le i \le n}E_{a_i}))$ is included in $\cup_{1\le 
i \le n}E_{a_i}$.  Deform $\beta$ as we did in the second 
step of the previous case. Hence we have an homotopy 
$\beta_t$, $t \in [0,v]$ for some $v<1$, such that 
$\beta_0=\beta$ and the number of self-intersection points of 
$\beta_v \cup (\cup_{1\le i \le n}E_{a_i})$ which are not 
double, is strictly smaller than  the number of 
self-intersection points of $\beta \cup (\cup_{1\le i \le 
n}E_{a_i})$ which are not double.  We can assume that all the 
intersection points of pairs of arcs of $\beta_t 
\cup(\cup_{1\le i \le n}E_{a_i})$ are transverse for each $t 
\in [0,v]$. Let us check now that for all $u \in [0,v]$, 
$\beta_u \cup (\cup_{1\le i \le n}E_{a_i})$ does not have 
bigons. Indeed, if for some $u \in [0,v]$, $\beta_u \cup 
(\cup_{1\le i \le n}E_{a_i})$ has a bigon $B$, we can follow 
$B$ in $\beta_t \cup (\cup_{1\le i \le n}E_{a_i})$. Since all 
intersection points of pair of arcs are transverse, there 
exists a bigon in $\beta_t \cup (\cup_{1\le i \le n}E_{a_i})$ 
for every $t \in [0,v]$. In particular, there exists a bigon 
in  $\beta \cup (\cup_{1\le i \le n}E_{a_i})$, a 
contradiction. (Note that the dotted regions of Figure 
\ref{beta}(b) are possible traces of bigons). Now we can 
extend the homotopy, removing at each step non-double 
self-intersection points. Thus the claim follows. 

By $(ii)$, $\VV_{\beta_1}$ is reduced. Thus, as above, one can prove that
$\VV_{\beta_1}=\cc(x_0x_1\dots x_{m-1})$.

Set $\alpha=\beta_1$. There exists $m$ oriented subarcs of $\alpha$, $B_0, B_1, \dots
B_{m-1}$ such that $\alpha=B_0B_1\dots B_{m-1}$ and for each $i$, the interior of $B_i$
does not intersect $\cup_{1\le i \le n}E_{a_i}$ and $B_i$ starts at $E_{x_i}$, runs on
the $\ovx_i$-side and ends at $E_{x_{i+1}}$, on the $x_{i+1}$-side. For each $i \in
\conj{0,1,\dots, m-1}$ there exists a unique arc $A_i$ joining two edges of $\Po$ such
that $\pi(A_i)=B_i$. Thus, $\alpha=\pi(A_0 A_1, \dots A_{m-1})$ satisfies the required
properties of a segmented representative without bigons.
\end{proof}

We will now show that there exists a segmented representative 
of a power of a primitive cyclic reduced word $\WW^r$, for 
which there exists $r-1$ self-intersection points which 
appear because the curve wraps around itself $r$ times. 
Moreover, we see that we can identify in which pair of 
segments of the segmented representative one can find these 
$r-1$ points and that the endpoints of every bigon are among 
these $r-1$ points.

\start{prop}{non primitive} Let $r \ge 1$ and let $\WW=\cc(y_0y_1\dots y_{m-1})^r$ be a
 reduced cyclic word  such that $\cc(y_0y_1\dots y_{m-1})$
 is primitive. Then there exists a segmented representative of $\WW$,
 $\alpha=\pi(A_0A_1\dots, A_{mr-1})$ such that
for every pair $i,j \in \conj{0,1,\dots, m r-1}$ the following are equivalent
\begin{numlist}
\item $A_i \cap A_j \ne \emptyset$   and $i \equiv j \pmod m$.
\item  $i=j$ or  $i=mr-1$ or $j=mr-1$.
\end{numlist}
Moreover, all the bigons of $\alpha$ have endpoints in the intersection of $A_{ m r-1}$
with $A_{m j-1}$ for some $j \in \conj{1,2,\dots, r-1}$.
\end{prop}
\begin{proof} By \propc{gon} there exists a segmented representative $\gamma$
of $\cc(y_0y_1\dots y_{m-1})$ which does not have bigons (see Figure \ref{arc
representatives}(a)). Assume that $\gamma=\pi(C_0C_1\dots C_{m-1})$. Let $A \subset \So$
be an annulus, having $\gamma$ as one of the boundary components. Subdivide this annulus
in $r-1$ parallel annuli. Let $\gamma_1, \gamma_2,\dots,\gamma_r$ be the boundary
components of these annuli. We can assume going along the arc $E_{y_0}$ in one of the two
possible directions, right after the initial arc of $\gamma_1$ (that is, the arc
``parallel'' to $C_0$ ),  one finds the initial arc of $\gamma_2$, then the initial arc
of $\gamma_2$  and so on. Each $\gamma_i$ is homotopic to $\gamma$, transversal to
$\cup_{1\le i \le n}E_{a_i}$ and such that for each $i \in \conj{1,2,\dots,r}$
$\gamma_i=\pi(D_{i,0}D_{i,1}\dots D_{i,m-1})$, where $D_{i,j}$ is an arc in $\Po$ from
the edge $\ovx_i$ to the edge $x_{i+1}$. Moreover, all the endpoints of the edges
$D_{i,s}$ are different for every $i \in \conj{0, 1, \dots, r-1}$ and $s \in
\conj{0,1,\dots, m-1}$   and if $i \ne j$, $D_{i,h} \cap D_{j, h}=\emptyset$.

\begin{figure}
\begin{pspicture}(14,4)
\rput(3.5,-0.3){$(a)$}\rput(10.5,-0.3){$(b)$}
 \rput{14}(3.5,2){\PstPolygon[unit=1.6,PolyNbSides=12]}
 \rput{14}(10.5,2){\PstPolygon[unit=1.6,PolyNbSides=12]}
\rput(-3.5,0){ \rput{0}(8.9,1.95){$a_1$} \rput(8,3.55){$\ova_1$}\rput(6,3.55){$a_2$}
\rput{0}(5.2,2){$a_3$}\rput{0}(6.05,0.55){$\ova_2$}\rput(8.1,0.55){$\ova_3$}}

\rput(3.5,0){ \rput{0}(8.9,1.95){$a_1$} \rput(8,3.55){$\ova_1$}\rput(6,3.55){$a_2$}
\rput{0}(5.2,2){$a_3$}\rput{0}(6.05,0.55){$\ova_2$}\rput(8.1,0.55){$\ova_3$}}

\rput(1,3){$\cc(a_1 \ova_2 \ova_3a_2)$}

\rput(7.5,3){$\cc((a_1 \ova_2 \ova_3a_2)^3)$}

\psline[linecolor=darkgray]{->}(4.2,3.4)(2.4,0.8)

\psline[linecolor=darkgray]{->}(2.6,3.28)(2,1.8)

\psline[linecolor=darkgray]{->}(4.4,0.7)(2.8,3.4)

\psline[linecolor=darkgray]{->}(2.8,0.6)(5,2.2)

\psline[linecolor=darkgray,linewidth=0.4pt]{->}(11.2,3.4)(9.4,0.83)
\psline[linecolor=darkgray,linewidth=0.4pt]{->}(11.3,3.3)(9.45,0.79)
\psline[linecolor=darkgray,linewidth=0.4pt]{->}(11.1,3.5)(9.35,0.89)

\psline[linecolor=darkgray,linewidth=0.4pt]{->}(9.5,3.2)(9,2.2)
\psline[linecolor=darkgray,linewidth=0.4pt]{->}(9.58,3.22)(9,2)
\psline[linecolor=darkgray,linewidth=0.4pt]{->}(9.68,3.28)(9,1.8)

\psline[linecolor=darkgray,linewidth=0.4pt]{->}(11.2,0.6)(9.74,3.35)
\psline[linecolor=darkgray,linewidth=0.4pt]{->}(11.4,0.7)(9.8,3.4)
\psline[linecolor=darkgray,linewidth=0.4pt]{->}(11.6,0.8)(9.86,3.45)

\psline[linecolor=darkgray,linewidth=0.4pt]{->}(9.7,0.7)(12.05,2.05)
\psline[linecolor=darkgray,linewidth=0.4pt]{->}(9.8,0.6)(12.05,2.35)
\psline[linecolor=darkgray,linewidth=0.4pt]{->}(9.9,0.5)(12.05,2.2)

\rput(13,3){$F_2$}\rput(13,2){$F_0$}\rput(13,1){$F_1$}

\psline[linestyle=dotted]{->}(12.8,3)(12,2.4)

\psline[linestyle=dotted]{->}(12.8,2)(12,2.25)

\psline[linestyle=dotted]{->}(12.8,1)(12,2)
\end{pspicture}
\caption{Minimal segmented representatives}\label{arc representatives}
\end{figure}
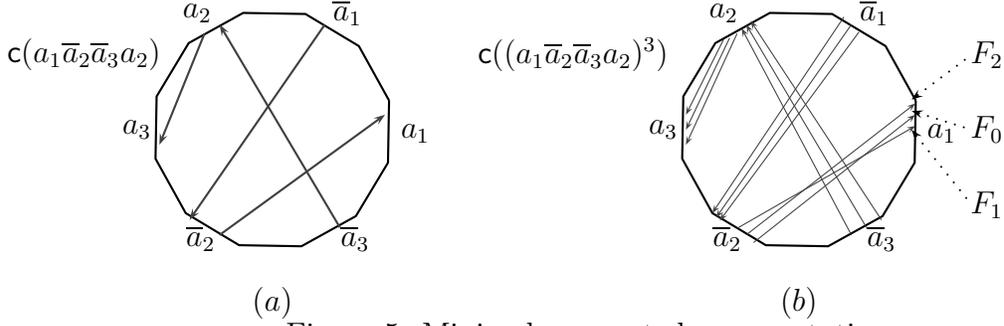

For each $j \in \conj{1,2,\dots,r-1}$, let $F_j$ be an arc in $\Po$ joining the initial
point of $D_{j,m-1}$ with the last point of $D_{j+1,m-1}$. We can assume that for every
$j, k \in \conj{1,2,\dots, r-1}$, if $j\ne k$ then $F_j \cap F_k=\emptyset$. Let $F_r$ be
an arc in $\Po$ joining the initial point of $D_{r,m-1}$ with the last point of
$D_{0,m-1}$.  We may also assume that each of the arcs $F_i$ intersects each of the arcs
$D_{h,l}$ in at most one point.

For each $i \in \conj{0,1,\dots, m r-1}$, write $i = m h +s$ where $0 \le s <m$. If $s
\ne m-1$, set $A_i=D_{h+1,s}$ and if $r=m-1$, $A_i=F_{h+1}$. By definition, $\alpha$ is a
loop. Moreover, $\alpha$ is a segmented representative.

Assume that  $\alpha=\pi(A_0A_1\dots A_{r m-1})$ has a bigon. Let $i, j$ and $k$ be as in
\lemc{arc bigon}. Write $i=m h_i + s_i$, $j= m h_j +s_j$, with $0\le s_i, s_j\le m-1$. We
claim that  $i$, $i+k$, $j$ and $j+k$ are congruent to $m-1$ modulo $m$.

Indeed, suppose that  \lemc{arc bigon}(1) holds. (If 
\lemc{arc bigon}(2) holds the proof can be completed with 
analogous arguments). Since $A_i \cap A_j \ne \emptyset$, 
then either $s_i=s_j=m-1$ or $s_i \ne s_j$. If $s_i \ne s_j$, 
$C_{s_i} \cap C_{s_j} \ne \emptyset$ and $C_{s_i+k} \cap 
C_{s_j+k} \ne \emptyset$. Also, $C_{s_i+h} \cap C_{s_j+h} 
=\emptyset$ if $h \in \conj{1,2,\dots, k-1}$. Then $\gamma$ 
has a bigon, a contradiction. So $s_i=s_j=m-1$. Since 
$A_{i+k} \cap A_{j+k} \ne \emptyset$ and $i+k$ and $j+k$ are 
congruent modulo $m$, $i+k$ and $j+k$ must be congruent to 
$m-1$ modulo $m$. 

By definition, if $A_{s m-1}$ intersects $A_{t m-1}$, then $s=r$ or $t=r$ and the proof
is complete.
\end{proof}

\subsection{Self-intersection points and linked pairs}\label{self-intersection}

We start to study by an example  the relation of pairs of subwords of a cyclic reduced
word  and self-intersection points of a segmented representative.

\start{ex}{motivation} Let $\OO$ and $\WW$ be as in 
\exc{first}. Since $a_2\ova_3$ and $\ova_3 a_1$ are subwords 
of $\WW$, any segmented representative of $\WW$ contains the 
projection of two transversal segments, $B_1$ from $\ova_2$ 
to $\ova_3$ and $B_2$, from $a_3$ to $a_1$  (see 
\figc{motive}(a)). One might guess that the occurrence of 
$a_2\ova_3$ and $\ova_3 a_1$ as subwords of a cyclic word 
will imply a self-intersection point in every representative 
of $\WW$. We will prove that a generalization of this holds, 
that is, certain pairs of subwords of two letters imply 
self-intersection points of the representatives.

\begin{figure}[ht]
\begin{pspicture}(14,4)

\rput{14}(2,2){\PstPolygon[unit=1.6,PolyNbSides=16]}
\rput{0}(0.3,2.7){$a_1$} \rput(1.4,3.68){$a_2$} \rput(2.7,3.65){$\ova_1$}
\rput(3.7,2.76){$\ova_2$}\rput(3.65,1.27){$a_3$}\rput(2.69,0.35){$a_4$}
\rput(1.32,0.35){$\ova_3$}\rput(0.4,1.2){$\ova_4$}
\psline[arrowinset=0,doublesep=0.0003pt,linecolor=darkgray]{->}(3.4,2.6)(1.4,0.6)
\psline[arrowinset=0,doublesep=0.0003pt,linecolor=darkgray]{->}(3.5,1.4)(0.6,2.6)

\rput{14}(6,2){\PstPolygon[unit=1.6,PolyNbSides=16]}
\rput{14}(9.18,2){\PstPolygon[unit=1.6,PolyNbSides=16]}
\rput{14}(12.3,2){\PstPolygon[unit=1.6,PolyNbSides=16]}
\rput(7.6,3.2){$a_1$} \pscurve[linestyle=dotted]{->}(7.6,3)(7,2.8)(7.5,2.2)
\rput(10.5,3.2){$a_1$} \pscurve[linestyle=dotted]{->}(10.7,3)(10.4,2.8)(10.6,2.2)
\rput(7.3,0.7){$\ova_1$}\rput(5.4,0.3){$a_3$}
\rput(8,0.4){$\ova_1$}\pscurve[linestyle=dotted]{->}(8,0.6)(8.2,0.8)(7.8,1.6)
\rput(10.3,0.2){$\ova_1$}\pscurve[linestyle=dotted]{->}(10.4,0.3)(11.2,1.3)(10.8,1.7)
\rput(11.1,3.3){$a_3$}\rput(11.2,0.6){$a_1$}

\pscurve[arrowinset=0,doublesep=0.0003pt,linecolor=darkgray]{->}(7.2,0.99)(7.2,1.4)(7.6,1.8)

\pscurve[arrowinset=0,doublesep=0.0003pt,linecolor=darkgray]{->}(5.4,0.6)(6.4,1.6)(7.6,2.1)

\psline[arrowinset=0,doublesep=0.0003pt,linecolor=darkgray]{->}(7.64,1.8)(10.7,1.8)

\psline[arrowinset=0,doublesep=0.0003pt,linecolor=darkgray]{->}(7.64,2.1)(10.7,2.1)

\pscurve[arrowinset=0,doublesep=0.0003pt,linecolor=darkgray]{->}(10.7,1.8)(11.2,2.6)(11.3,3.2)

\pscurve[arrowinset=0,doublesep=0.0003pt,linecolor=darkgray]{->}(10.7,2.1)(11.3,1.6)(11.3,0.8)

\rput(2,-0.1){$(a)$}\rput(9,-0.1){$(b)$}

\end{pspicture}\caption{$\exc{motivation}$}\label{motive}
\end{figure}

Now, consider the pair of subwords of $\WW$, $a_1a_1$ and 
$\ova_3a_1$, see Figure \ref{motive}(b). Since both the 
segments corresponding to this pair of subwords  land in the 
edge $a_1$, the occurrence of these two subwords does not 
provide enough information to deduce the existence of a 
self-intersection point.  In order to understand better this 
configuration of segments, we prolong the subwords starting 
with $a_1a_1$ and $\ova_3a_1$ until they have different 
letters at the beginning and at the end. Then we study how 
the arcs corresponding to these subwords intersect. So for 
instance, in our example we get $\ova_3 a_1 a_1$ and $a_1 a_1 
a_3$, and we  see a self-intersection point (Figure 
\ref{motive}(b)). We will see that certain pairs of subwords 
of $\WW$ imply self-intersection points of representatives. 
\end{ex}

Let $\alpha=\pi(A_0A_1\dots,A_{m-1})$ be a segmented 
representative of a reduced cyclic word $\WW=\cc(x_0x_1\dots 
x_{m-1})$. By an {\em arc} of  $\alpha$  we mean a finite 
subsequence of segments of the infinite sequence $\prod_{i 
\in \N}\pi(A_1)\pi(A_{2})\dots \pi(A_{m})$.  The {\em 
underlying subword of the arc $\pi(A_iA_{i+1}\dots A_{i+j})$} 
is the subword $x_i x_{i+1}\dots x_{i+j+1}$ of 
$\WW^{\infty}$.

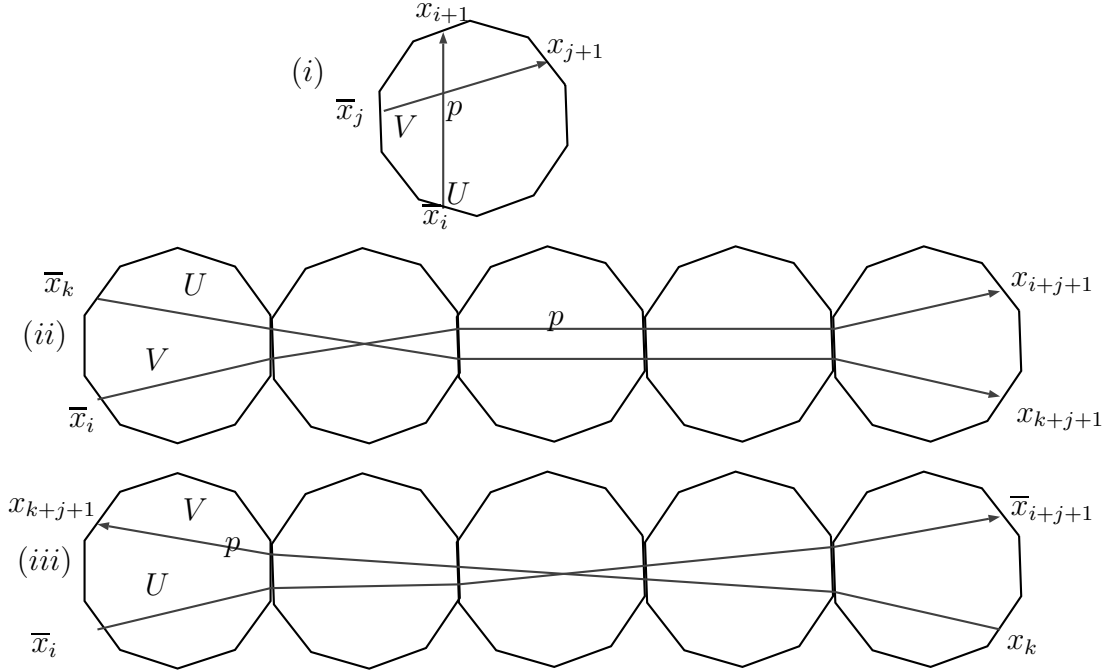
\begin{figure}[ht]
\begin{pspicture}(14,9)
\rput(6,6.4){$U$}\rput(5.3,7.3){$V$} \rput(2.5,5.2){$U$} \rput(2,4.2){$V$} \rput(0,-3){
\rput(2.5,5.2){$V$} \rput(2,4.2){$U$} }
\rput(0,3){\rput(1,0.4){$\ovx_{i}$} \rput(0.7,2.2){$\ovx_{k}$}}
\rput(0,3){\rput(14,0.4){$x_{k+j+1}$} \rput(13.9,2.2){$x_{i+j+1}$}}
\rput(0.5,0.4){$\ovx_{i}$} \rput(0.6,2.2){$x_{k+j+1}$} \rput(13.5,0.4){$x_{k}$}
\rput(13.9,2.2){$\ovx_{i+j+1}$}
\rput(4,8){$(i)$}\rput(0.5,4.5){$(ii)$}\rput(0.5,1.5){$(iii)$} \rput(-1,0){
\rput(1,6){\psline[arrowinset=0,linecolor=darkgray]{->}(5.8,0.2)(5.8,2.55)
\psline[arrowinset=0,linecolor=darkgray]{->}(5.01,1.5)(7.2,2.156)
\rput(5.7,0.1){$\overline{x}_i$} \rput(5.8,2.8){${x}_{i+1}$}
\rput(4.57,1.47){$\overline{x}_j$} \rput(7.56,2.3){$x_{j+1}$}

\rput{20}(6.2,1.4){\PstPolygon[unit=1.3,PolyNbSides=10]}}
\rput(2,0){ \rput(0,0){\sequence} \rput{0}(0,3){\sequence}
\psline[arrowinset=0,linecolor=darkgray]
{->}(0.2,3.66)(2.5,4.2)(5,4.6)(7.5,4.6)(10,4.6)(12.2,5.1)
\psline[arrowinset=0,linecolor=darkgray,doublesep=0.0003cm]
{->}(0.2,5)(2.5,4.6)(5,4.2)(10,4.2)(12.2,3.7)
\psline[arrowinset=0,linecolor=darkgray,doublesep=0.0003cm]
{->}(0.2,0.6)(2.5,1.15)(5,1.2)(10,1.7)(12.2,2.1)
\psline[arrowinset=0,linecolor=darkgray,doublesep=0.0003cm]
{<-}(0.2,2)(2.5,1.6)(10,1.1)(12.2,0.6)}}
\rput(5.95,7.5){$p$}\rput(7.3,4.7){$p$}\rput(3,1.7){$p$}
\end{pspicture}
\caption{Strips and semiparallel pairs.}\label{strip}
\end{figure}

\start{defi}{semi} Let $U$, $V$  be a pair of arcs of 
segmented representatives of $\alpha$ and $\beta$ 
respectively,  such that the underlying words of $U$ and $V$ 
are $x_i x_{i+1}\dots x_{i+j+1}$, and $y_k y_{k+1}\dots 
y_{k+j+1}$, respectively. Moreover, assume that exactly one 
of the following holds (see Figure \ref{strip}). 
\begin{numlist}
\item  $x_{i}\ne y_j$, $x_{i+j+1}\ne y_{k+j+1}$ and
$x_{i+1}x_{i+2}\dots x_{i+j}=y_{k+1}y_{k+2}\dots y_{k+j}$.
\item $x_{i}\ne \ovy_{k+j+1}$, $x_{i+j+1}\ne \ovy_{k}$ and  $x_{i+1}x_{i+2}\dots x_{i+j}=
\ovy_{k+j}\ovy_{k+j-1}\dots \ovy_{k+1}$.
\end{numlist}
If $\alpha \ne \beta$  (resp. $\alpha=\beta$) we say that 
$\conj{U, V}$ is a {\em semiparallel pair of arcs of $\alpha$ 
and $\beta$ (resp. of $\alpha$)} parallel in case (1) and 
antiparallel in case (2). 
\end{defi}

Let $\VV$ be a  primitive reduced cyclic word of length $m$, 
let $r \ge 1$ let $\WW=\VV^r$. We call a segmented 
representative of $\WW$ as in \propc{gon} if $r=1$ and as in 
\propc{non primitive} if $r \ge 1$  {\em minimal}. Let $r>1$ 
and $\alpha$ be a minimal segmented representative of $\WW$. 
We denote by $\ap_\alpha$ the set of intersection points of 
$\pi(A_{rm-1})$ with $\pi(A_{km-1})$, for $k \in 
\conj{1,2,\dots, r-1}$ and by $\ai_\alpha$ the set of 
self-intersection points of $\alpha$ not in $\ap_\alpha$. 
Hence the set of self-intersection points of a minimal 
segmented representative $\alpha$ is the disjoint union of 
$\ap_\alpha$ and $\ai_\alpha$. When $r=1$, $\ap_\alpha$ is 
empty by definition and $\ai_\alpha$ is the set of 
self-intersection points of $\alpha$.

\start{lem}{pair} Let $\WW=\cc(x_0x_1\dots x_{m-1})$ be a 
cyclically reduced cyclic word and let $\alpha$ be a minimal 
segmented representative of $\WW$. Let $p$ be a 
self-intersection point of $\alpha$, such that $p \in 
\ai_\alpha$. Then there exists a unique semiparallel pair of 
arcs of $\alpha$, $\conj{U,V}$, such that $p \in U \cap V$. 
Moreover, if $U=\pi(A_iA_{i+1}\dots A_{i+j})$ and  
$V=\pi(A_kA_{k+1}\dots A_{k+j})$ then $0 \le j \le 
\len(\WW)-1$ and there exists a unique $u \in 
\conj{0,1,\dots, j}$ such that $p \in \pi(A_{i+u})$ and $p 
\in \pi(A_{j+u})$ in parallel case and  $p \in  
\pi(A_{k+j-u})$) in the antiparallel case. 
\end{lem}

\begin{proof} Suppose that $p \in \pi(A_r)\cap \pi(A_s)$. Since $\WW$ is cyclically reduced exactly one
of the following holds.
\begin{romlist}
\item $\{\ox_r,x_{r+1}\} \cap\{\ox_s,x_{s+1}\}=\emptyset$
\item $\{\ox_r,x_{r+1}\} \cap\{\ox_s,x_{s+1}\}\ne \emptyset$, and $\ox_r=\ox_s$ or $x_{r+1}=x_{s+1}$.
\item $\{\ox_r,x_{r+1}\} \cap\{\ox_s,x_{s+1}\}\ne \emptyset$, and $\ox_r=x_{s+1}$ or $x_{r+1}=\ox_s$.
\end{romlist}

If $(i)$ holds, setting $\conj{U,V}=\conj{\pi(A_r), \pi(A_s)}$,taking  $i=r$, $k=s$ and
$j=0$, the result holds.

In the cases $(ii)$ and $(iii)$ one keeps adding segments before and after $A_r$ and
$A_s$, until finding arcs landing in different edges at the beginning and at the end.
More precisely, assume that $(ii)$ holds. By \lemc{well} and since $p \notin \ap_\alpha$,
there exist integers $t,l\ge -1$, such that $1 \le t+l < \len(\WW)-1$ and $t+l$ is
maximum with the property that $$x_{r-t}x_{r-t+1}\dots x_{r+l}=x_{s-t}x_{s-t+1}\dots
x_{s+l}.$$

Here, we set $\conj{U,V}=\conj{\pi(A_{r-t-1}A_{r-t}\dots
A_{r+l}),\pi(A_{s-t-1}A_{s-t}\dots A_{s+l})}$. Clearly, 
\defc{semi}(1) holds. Thus by taking $i=r-t-1 \pmod m$, 
$k=s-t-1 \pmod m$, and  $u=t+1$, $j=t+l+1$,   the proof is 
complete for this case. 

Finally, assume that $(iii)$ holds. Let $t,l$ be  
non-negative integers such that $t+l\ge 1$ and $t+l$ is the 
maximum positive integer  such that $$x_{r-t+1}x_{r-t+2}\dots 
x_r x_{r+1}\dots x_{r+l}=\ox_{s+t}\ox_{s+t-1}\dots 
\ox_{s+1}\ovx_s \dots \ox_{s-l+1}.$$ Here, we define  
$U=\pi(A_{r-t}A_{r-t+1}\dots A_{r+l})$ and 
$V=\pi(A_{s-l}A_{s-l+1}\dots A_{s+t})$. Hence $\conj{U,V}$ 
satisfies \defc{semi}(2), and so taking $i=r-t-1 \pmod m$, 
$k=s-l \pmod m$, and $u=t$ and $j=t+l+1$ the proof of the 
lemma is complete. 
\end{proof}

We associate to linear word $x_i x_{i+1}\dots x_{i+j}$ a 
sequence of copies of $\Po$ glued in a certain  way:  
Consider  $j-1$ copies $\Po^1, \Po^2, \dots, \Po^{j-1}$, 
 of the polygon $\Po$. Set $\S=(\Po^1 \cup \Po^2 \cup \dots 
\Po^{j-1})/\sim$, where $\sim$ is the equivalence relation 
generated by the pairs $(y,z)$ for which  there exists $h \in 
\conj{1,2,\dots,j-1}$ such that  $y$ is in the edge of 
$\Po^{i+h-1}$ labeled with $x_{i+h}$, $z \in \Po^{i+h}$, and 
$z$ is in the edge of $\Po^{i+h}$ labeled with $\ovx_{i+h}$ 
and $\pi(y)=\pi(z)$. Such an $\S$ is a {\em strip of the word 
$x_i x_{i+1}\dots x_{i+j}$} (see Figure \ref{strip}). Observe 
that for each $h \in \conj{1,2,\dots,j}$, $\Po^h$ is embedded 
in $\S$. Thus we can think of $\Po^h$ as a subset of $\S$.  
The map $\flecha{\Po^1 \cup \Po^2 \cup \dots \Po^j}{}{\So}$ 
which restricted  to each of the copies of $\Po$ is the 
projection $\pi$, induces a map, $\map{\Pi}{\S}{\So}$.

Let $C(\{U,V\})$ denote the subset of points of $\ai_\alpha$ 
which are associated to $\{U,V\}$ by \lemc{pair}.

Denote the image of a map $\beta$ by $\im(\beta)$. By 
\lemc{pair} we have:

\start{lem}{lift} Let $U=\pi(A_iA_{i+1}\dots A_{i+j})$, 
$V=\pi(A_kA_{k+1}\dots A_{k+j})$
be a semiparallel pair of arcs of a segmented 
representative $\alpha$ and let $\S$ be the
strip of the underlying word of $U$.
 Then there exist two continuous maps, $\map{\mu,
\nu}{[0,1]}{\S}$ such that
\begin{numlist}
\item $\Pi(\im(\mu))=U$ and $\Pi(\im(\nu))=V$. Moreover, 
for each $h \in \conj{0,1,2,\dots,j}$,
$\Pi(\im(\mu) \cap \Po^h )=\pi(A_{i+h})$ and if the pair of 
arcs $\conj{U,V}$ are parallel (resp. antiparallel) then 
$\Pi(\im(\nu) \cap \Po^h )=\pi(A_{k+h})$ (resp. $\Pi(\im(\nu) 
\cap \Po^h )=\pi(A_{k+j-h})$). 

\item $\Pi(\im(\mu) \cap \im(\nu))=C(\{U,V\})$. Furthermore, the 
restriction of $\Pi$ to $\im(\mu) \cap \im(\nu)$ is a 
bijection between $\im(\mu) \cap \im(\nu)$ and $C(\{U,V\})$. 
\end{numlist}
\end{lem}

\start{theo}{bigons} Let $\WW=\cc(x_0x_1\dots x_{m-1})$ be a 
cyclically reduced cyclic word in the letters of $\an$    and 
let $\alpha$ be a minimal segmented representative of $\WW$. 
Then self-intersection points of $\alpha$ in $\ai_\alpha$ 
correspond bijectively to the set of pairs of linked pairs of 
$\WW$ of the form $\conj{(\ap, \aq), (\aq,\ap)}$. 
\end{theo}

\begin{proof} We claim that  different points in $\ai_\alpha$ cannot be assigned to the same pair
of arcs by \lemc{pair}. Indeed let $U=\pi(A_iA_{i+1}\dots 
A_{i+j})$, $V=\pi(A_kA_{k+1}\dots A_{k+j})$ be a semiparallel 
pair of arcs and assume that $\conj{U,V}$ are assigned to 
more than one point. Take two of these points $p, q$. By 
\lemc{pair} we can assume that there exist $s, t \in 
\conj{0,1,\dots,j}$ such that $s<t$, $p \in \pi(A_{i+s})$, 
$q\in \pi(A_{i+t})$ and there are no points assigned to 
$\conj{U,V}$ in the arc $\pi(A_{i+s+1}\dots A_{i+t-1})$. Let 
$\mu$, $\nu$ and $\S$ be as in \lemc{lift}. Let $P, Q \in \S$ 
be such that $\Pi(P)=p$ and $\Pi(Q)=q$.  The subarc of 
$\im(\mu)$ from $P$ to $Q$ and the subarc of $\im(\nu)$ from 
$Q$ to $P$ bound a disk. Therefore the subarc of $U$ from $p$ 
to $q$ and the subarc of $V$ from $q$ to $p$ bound the image 
of a disk. Since the bigons of $\alpha$  have endpoints in 
$\ap_\alpha$, $p, q \notin \ai_\alpha$, a contradiction. So 
the claim holds. 

Now, let $p \in \ai_\alpha$, let $\conj{U,V}$ be the pair of 
arcs assigned to $p$ by \lemc{pair} and let $\ap$ and $\aq$ 
be the underlying words of $U$ and $V$ respectively. Let 
$\mu$, $\nu$ and $\S$ be as in \lemc{lift}. By the above 
claim $p$ is the only point in $\im(\mu) \cap \im(\nu)$.
Therefore $(\ap, \aq)$ is a linked pair (see Figure 
\ref{strip}).

Conversely, if $(\ap, \aq)$ is a linked pair of occurrences 
of subwords of $\WW$, $U$ and $V$ are the arcs with 
underlying words $\ap$, $\aq$ respectively, then $U$, $V$ is 
a semiparallel pair of arcs. Hence we can apply \lemc{lift}. 
By the definition of linked pair, one has $\im(\mu) \cap 
\im(\nu)$ contains a single point $P$, (see Figure 
\ref{strip}). The pair of arcs assigned to $P$ by \lemc{pair} 
is $\conj{U,V}$. Hence the proof of the theorem is complete. 

\end{proof}

\start{rem}{CL1} By  construction,  a minimal segmented 
representative of a primitive cyclic reduced word $\WW$ has 
the minimal number of self-intersection points in its free 
homotopy class. Hence by \theoc{bigons}, the minimal number 
of self-intersection points of pairs of representatives $\WW$ 
equals half of the cardinality of the set of linked pairs of 
$\WW$, $\bb_1(\WW)$. This is an equivalent form of a result 
found in stages by Birman and Series \cite{BS}, and Cohen and 
Lustig in \cite{CL}, (see also \cite{Tan}). 
\end{rem}

\subsection{Study of intersection points of a pair of segmented representatives}

A pair of representatives $\alpha, \beta$ of a pair of reduced cyclic words $\VV$ and
$\WW$ such that all intersections are transverse double points has {\em minimal
intersection} if every pair of representatives of $\VV$ and $\WW$, such that their
intersection consists in finitely many transverse double points,  has at least as many
intersection points as $\alpha$ with $\beta$.

 The next result is analogous to  \propc{non primitive}.

\start{prop}{good pair} For each pair of reduced cyclic words 
$\VV$ and $\WW$ there exists a pair of segmented 
representatives $\alpha$ and $\beta$  with the following 
properties: The union of $\alpha$ and  $\beta$  does not have 
bigons, the endpoints of the arcs of $\alpha$ do not 
intersect the endpoints of the arcs of $\beta$ and there are 
no triple intersection points between $\alpha$ and $\beta$. 
Furthermore, if $\VV$ and $\WW$ are not powers (positive or 
negative) of the same cyclic word then $\alpha$ and $\beta$ 
have minimal intersection. 

Finally, if $\VV$ and $\WW$ are (positive or negative) powers  of a primitive cyclic word
of length $k$, $\alpha=\pi(A_0A_1\dots A_m)$ and $\beta=\pi(B_0B_1\dots B_l)$, and if $k$
divides $i-j$, for some $i \in \conj{0,1,\dots,m}$, and $j \in \conj{0,1,\dots,l}$ then
$\pi(A_i) \cap \pi(B_j)=0$.
\end{prop}
\begin{proof} First assume that $\VV=\VV_1^r$, $\WW=\WW_1^s$ where $\VV_1$ and $\WW_1$
are two different primitive cyclic words,  $\VV_1\ne \ov{\WW}_1$ and $r, s \ge 1$.

Construct two representatives of $\VV_1$ and $\WW_1$ 
$\alpha_1$ and $\beta_1$, as in \propc{gon}, considering two 
geodesic represented by $\VV_1$ and $\WW_1$, respectively. 
Observe that $\alpha_1 \cup \beta_1$ does not have bigons. If 
there are triple intersection points of $\cup_{1\le i \le 
n}E_{a_i} \cup \alpha_1 \cup \beta_1$, deform $\alpha \cup 
\beta$ as in \propc{gon} to  $\alpha_1 \cup \beta_1$ where 
$\alpha_1$ and $\beta_1$ are a pair of curves such that 
$\alpha_1 \cup \beta_1$ does not have bigons, nor triple 
intersection points. Now, construct $\alpha$ and $\beta$ 
using $\alpha_2$ and $\beta_2$ respectively as in the proof 
of \propc{non primitive}. Observe that if $\alpha \cup \beta$ 
has a bigon then $\alpha_2\cup\beta_2$ has a bigon.

The minimal number of intersection points of $\alpha$ and 
$\beta$, equals the number of intersection points of 
$\alpha_1^r$ and $\beta_1^s$ (where  the  $p$-th power of a 
representative means run along the representative $p$ times). 
Using this, it is not hard to see that $\alpha$ and $\beta$ 
have minimal intersection. 

Now, assume that $\VV$ and $\WW$ are powers (positive or 
negative) of a word $\XX$. We can assume that $\XX$ is 
primitive. Let $\gamma$ be a segmented representative of 
$\XX$. Consider an  annulus $N$ in $\So$  around $\gamma$, so 
that both boundary components if $N$ are also segmented 
representatives of $\XX$. Then $N \setminus \gamma$ is the 
disjoint union of two annuli, $N_1$ and $N_2$. Let $\alpha$ 
and $\beta$ be arcs representatives of $\VV$ and $\WW$ such 
that $\alpha \subset N_1$ and $\beta \subset N_2$, and every 
arc of $\alpha$ intersects every arc of $\beta$ in at most 
one point. Notice that if $\alpha$ and $\beta$ have a bigon, 
then $\gamma$ has a bigon. So the proof is complete. 
\end{proof}

Let $\VV$ and $\WW$ be a pair of reduced cyclic words. A pair 
of representatives $\alpha$ and $\beta$ as in \propc{good 
pair} will be called \emph{good representatives of $\VV$ and 
$\WW$} 

The proof of the next result is  analogous to the proof  of 
\lemc{pair} but one needs to consider the following two cases 
separately: either the two curves are powers (positive or 
negative) of the same primitive word or they are not. In the 
former case, one also needs to apply \propc{finite length}. 

\start{lem}{pair two} Let $\VV=\cc(y_0y_1\dots y_{h-1})$ and 
$\WW=\cc(x_0x_1\dots x_{m-1})$ be  cyclically reduced cyclic 
words and let $\alpha$, $\beta$ be  good representatives of 
$\VV$ and $\WW$. Let $p\in \alpha \cap \beta$. Then there 
exists a unique semiparallel pair of arcs of $\alpha$ and 
$\beta$, $\conj{U,V}$, such that $p \in U \cap V$. Moreover, 
if $U=\pi(A_iA_{i+1}\dots A_{i+j})$ and $V=\pi(B_k 
B_{k+1}\dots B_{k+j})$ then  there exists a unique $u \in 
\conj{0,1,\dots, j}$ such that $p \in \pi(A_{i+u})$ and $p 
\in \pi(B_{j+u})$ if $U$ and $V$ are parallel and  $p \in  
\pi(B_{k+j-u})$), if $U$ and $V$ are antiparallel. 

\end{lem}

The next result is the equivalent of \theoc{bigons} for pairs of cyclic words, and the
proof uses analogous ideas.

\start{theo}{bigons for a pair} Let $\VV$ and $\WW$ be two 
primitive reduced cyclic words. Let $\alpha$ and $\beta$ be a 
good pair of representatives of $\VV$ and $\WW$ (as in 
\propc{good pair}). Then there is a one to one correspondence 
between intersection points of $\alpha$ and $\beta$ and 
$\bb_2(\VV,\WW)$, the set of linked pairs of $\VV$ and $\WW$. 
\end{theo}

The next remark is analogous to  \remc{CL1}.

\start{rem}{CL2} Let $\VV$ and $\WW$ be primitive reduced cyclic words.  The minimal arc
representatives of $\VV$ and $\WW$ have minimal  intersection, and then, the minimal
number of intersection points  of representatives of $\VV$ and $\WW$ equals the
cardinality of the set of linked pairs of $\VV$ and $\WW$, $\bb_2(\VV,\WW)$. An
equivalent form of this result was obtained by Cohen and Lustig, \cite{CL}.
\end{rem}

\section{The isomorphism}\label{GT}

Let $\bri{\,, }$ and $\cib$ be the bracket and cobracket of the Lie bialgebra structure
defined geometrically in Appendix \ref{definition of the Lie bialgebra} on the vector
space generated by non-trivial free homotopy classes of curves on an oriented smooth
surface (Appendix \ref{curves}). By \lemc{correspondence}, $\bri{\,,}$ and $\cib$ are
defined on $\V$, the vector space generated by all cyclically reduced cyclic words on the
alphabet $\an$. In other words, $(\V, \bri{\,,},\cib)$ is a Lie bialgebra. On the other
hand, we have defined in Section~\ref{words} two linear maps $\flecha{\V\otimes
\V}{\bra{,}}{\V}$ and $\flecha{\V}{\cob{ }}{\V\otimes \V}$. We will show that the
geometrically defined maps $\bri{\,,}$ and $\cib$ coincide with the combinatorially
defined maps $\bra{\,,}$ and $\cob$.

\start{prop}{isomorphism} For each $\VV, \WW \in \V$, 
$\cob(\VV)=\cib(\VV)$ and $\bra{\VV,\WW}=\bri{\VV,\WW}$. 
\end{prop}
\begin{proof} We prove that $\cob(\VV)=\cib(\VV)$. Write $\VV=\WW^r$
where $\WW$ is a  primitive reduced cyclic word. Let $\alpha$ be a segmented
representative of $\VV$ and let $S$ be the   set of self-intersection points of $\alpha$.
Recall that  $S$ is the disjoint union of two subsets, $\ai_\alpha$, and $\ap_\alpha$
where $\ap_\alpha$ is the set of intersection points of $A_{rm}$ with $A_{km}$, for $k
\in \conj{1,2,\dots, r-1}$. By \theoc{bigons}, $\ai_\alpha$ is in correspondence with the
pairs of linked pairs of $\WW$ of the form $\conj{(\ap,\aq), (\aq,\ap)}$. For each $q \in
S$, if $\VV^q_1$ and $\VV^q_2$ are as in Appendix \ref{curves} then $ \cib(\VV)= \sum_{q
\in S} \VV^q_1 \otimes \VV^q_2-\VV^q_2 \otimes \VV^q_1.$  By definition of $\ap_\alpha$,
$\sum_{q\in \ap_\alpha} \VV^q_1 \otimes \VV^q_2-\VV^q_2 \otimes\VV^q_1 =
\sum_{i=1}^{r-1}\WW^i\otimes\WW^{r-i}-\WW^{r-i}\otimes\WW^{i}=0.$

Consider $p \in \ai_\alpha$ and let $(\ap, \aq) \in 
\bb_1(\VV)$ be the linked pair assigned to $p$ by 
\theoc{bigons}. If $\sign(\ap, \aq)=1$ then 
$\WW_1^p=\cob_1(\ap,\aq)$ and $\WW_2^p=\cob_2(\ap,\aq)$ and 
if  $\sign(\ap, \aq)=-1$ then $\WW_1^p=\cob_2(\ap,\aq)$ and 
$\WW_2^p=\cob_1(\ap,\aq)$. Hence $\cib(\VV)=\sum_{q \in 
S_\OO}\WW^q_1\otimes \WW^q_2-\WW^q_2 \otimes 
\WW^q_1=\cob(\VV).$

The proof of $\bra{\VV,\WW}=\bri{\VV,\WW}$ is analogous.
\end{proof}

By Propositions~\ref{isomorphism} and \refc{involutive}, we have the following.

\start{theo}{Main} $(\V,\bra{\,,},\cob)$ is an involutive Lie bialgebra.
\end{theo}

\section{Applications}\label{applications}

A free homotopy class is {\em simple} if it has a simple 
representative, i.e., a non-selfintersecting representative. 
As simple closed curves, simple free homotopy classes satisfy 
statements that do not hold for the rest of the free homotopy 
classes. The extension of Goldman's result, \theoc{no 
cancellation}, is one of these. On the other hand, by 
definition, if a free homotopy class is simple then its 
cobracket is zero. Using the algorithm we describe in this 
paper, we developed a program in $C^{++}$ that, given the 
surface symbol $\OO$,  computes the bracket and  the 
cobracket of reduced cyclic words. By running this program we 
found classes with cobracket zero which are not powers of a 
simple class, some of which we list in Examples \ref{counter} 
and \refc{more} . By means of the same programs, we found 
examples of pairs of non-simple classes with bracket zero, 
and no disjoint representatives (\exc{and more}). 

\subsection{Topological proof of Goldman's result}\label{answer}

The next lemma will be used in the proof of \theoc{no cancellation}.

\start{lem}{hurra} Let $\VV$ be a cyclically reduced cyclic 
word and let $x \in \an$. If $(\ap_1, \aq_1)$, $(\ap_2,\aq_2) 
\in \bb_2(\VV,\cc(x))$ are linked pairs and $\br(\ap_1, 
\aq_1)=\br(\ap_2,\aq_2)$ then $(\ap_1,\aq_1)=(\ap_2,\aq_2)$ 
and so $\sign(\ap_1,\aq_1)=\sign(\ap_2,\aq_2)$. \end{lem} 

\begin{proof}

Linked pairs of $\VV$ and $\cc(x)$ have the form $(u_1x^ku_2,x^{k+2})$  and
$(u_1\ovx^ku_2,x^{k+2})$  where $u_1x^ku_2$ and $u_1\ovx^ku_2$ respectively are subwords
of $\VV$, $u_1\ne x$, $u_2 \ne x$ and $k \ge 0$ (in the first case,
\defc{def}$(1)$ or $(2)$ holds and in the second case, \defc{def}$(3)$ holds).

Observe that if   $(u_1\ovx^ku_2,x^{k+2})$ (resp. $(u_1x^ku_2,x^{k+2})$)) is a linked
pair then $$\len(\br(u_1\ovx^ku_2,x^{k+2}))=\len(\VV)+1 \mbox{ (resp.
$\len(\br(u_1x^ku_2,x^{k+2}))=\len(\VV)-1)$}.$$

Hence if $(\ap_1, \aq_1),$ $(\ap_2,\aq_2) \in 
\bb_2(\VV,\cc(x))$ are linked pairs such that $\br(\ap_1, 
\aq_1)=\br(\ap_2,\aq_2)$ then either $\ap_1=y_1x^ky_2$ and 
$\ap_2=z_1x^jz_2$ or $\ap_1=y_1\ovx^ky_2$ and 
$\ap_2=z_1\ovx^jz_2$. 

Assume that $\ap_1=y_1\ovx^ky_2$, $\ap_2=z_1\ovx^jz_2$,
$\br(\ap_1,x^{k+2})=\br(\ap_2,x^{j+2})$, where $k, j \ge 1$. 
There exist (possibly empty) linear words $\XXX$, $\YYY$ such 
that $\WW=\cc(\ap_1\XXX\ap_2\YYY)$, where, if $\XXX$ (resp. 
$\YYY$) is empty then  the last letter of $\ap_1$ (resp. 
$\ap_2$) can coincide with the first letter of $\ap_2$ (resp. 
$\ap_1$). Hence $\br(\ap_1,x^{k+2})=\cc(y_1x^{k-1}y_2\XXX 
\ap_2\YYY)\mbox{ and } \br(\ap_2,x^{j+2})=\cc(\ap_1\XXX z_1 
x^{j-1}z_2\YYY).$

If $\WW$ is a cyclic word, denote by $\CC_{\WW}$ the set of 
subwords of $\WW$ of the form $u_1\ovx^h u_2$ with $u_1, u_2 
\ne x$ and $h \ge 0$. Thus, $$\CC_\WW \setminus\conj{\ap_1, 
\ap_2}=\CC_{\br(\ap_1,x^{k+2})} \setminus\conj{y_1x^{k-1}y_2, 
\ap_2}=\CC_{\br(\ap_2,x^{j+2})} \setminus\conj{\ap_1, 
z_1x^{j-1}z_2}.$$ Since 
$\br(\ap_1,x^{k+2})=\br(\ap_2,x^{j+2})$, we have that 
$\CC_{\br(\ap_1,x^{k+2})}=\CC _{\br(\ap_2,x^{j+2})}.$ Hence 
$$\conj{y_1x^{k-1}y_2, \ap_2}=\conj{\ap_1, z_1x^{j-1}z_2}$$ 
and then, finally, $\ap_1=\ap_2$, as desired.

In order to prove the above result for the other case, one considers $\CC_{\WW}$, the set
of subwords of $\WW$ of the form 
$u_1x^h u_2$ with $u_1, u_2 \ne x$ and $h \ge 0$ and
proceeds as above.
\end{proof}

\start{lem}{correction} Let $\Sigma$ be an oriented  surface with  boundary and let $\lambda$ be a simple closed curve non-homologous to zero. Then there exists an alphabet $\an$, a surface word $\OO$ for that alphabet and a homeomorphism  $\map{\rho}{\Sigma}{\Sigma_\OO}$  such that the cyclic word for the free homotopy class of $\rho(\lambda)$ has one letter.

\end{lem} 

\begin{proof} One chooses a non-selfintersecting arc $a_1$ transversal to $\lambda$, with endpoints on boundary components and intersecting $\lambda$ exactly once. Then one cuts $\Sigma$ open along $\lambda$ and $a_1$.  
Now, one studies the separating and non-separating case.

In the separating case, this procedure yields  two surfaces. Then one continues cutting open these two surfaces along non-selfintersecting arcs $a_2, a_3, \dots a_n$ with endpoints in the boundary components of $\Sigma$,  until obtaining two disks.

In the non-separating case one cuts the surface open along non-selfintersecting arcs $a_2, a_3, \dots a_{n-1}$ with endpoints in the boundary components of $\Sigma$, until obtaining a disk. After gluing this disk along $\lambda$, one obtains a cylinder. Then  one chooses an   arc $a_n$ with endpoints  in the boundary components of $\Sigma$ such that cutting the cylinder open along $a_n$ yields a disk.

Denote by $\Po$ the disk obtained by cutting $\Sigma$ along $a_1, a_2, \dots a_n$.
Then taking the identity map as the homeomorphism $\rho$ completes the proof.

\end{proof}

 We state the next theorem  using the
correspondence between free homotopy classes and reduced cyclic words given by
\lemc{correspondence}.

 \start{theo}{no cancellation} Let $\VV$ and $\WW$ be cyclic reduced words and  such
that $\VV$ has a simple representative which is non-homologous to zero. Then there exists two  representatives $\alpha$
and $\beta$ of $\VV$ and $\WW$ respectively such that the bracket of $\VV$ and $\WW$
computed using the intersection points of  $\alpha$ and $\beta$ does not have
cancellation. In other words, the number of terms (counted with multiplicity) of
$\bri{\VV, \WW}$ equals the minimal number of intersection points of representatives of
$\VV$ and $\WW$.
\end{theo}
\begin{proof} 
We can assume that $\VV$ and $\WW$ are not powers of the same word because if  $\WW$ is a power of $\VV$, then $\VV$ and $\WW$ have disjoint representatives. 

By \propc{good pair}, $\VV$ and $\WW$ have minimal intersection. Moreover, the number of terms of the brackets is the minimal number of intersection points. Therefore, by \propc{isomorphism} if the bracket of two words in a certain alphabet has cancelation, then the bracket of the corresponding words in another alphabet also has cancelation.

Thus by \lemc{correction} we can assume that $\VV$ has only one letter. Then the result follows from \lemc{hurra}.

\end{proof}

The next is a particular case of \theoc{no cancellation}.

\start{cory}{cory} Let $\VV$ and $\WW$ be two free homotopy 
classes of curves on a surface with boundary and assume that 
$\VV$ has a simple representative which is non-homologous to zero. Then the bracket of $\VV$ 
and $\WW$ is zero if and only if  $\VV$ and $\WW$ have 
disjoint representatives. 
\end{cory}

The next example shows that the hypothesis that one of the classes is simple cannot be
omitted from \theoc{no cancellation}.

\start{ex}{and more} Let $\OO=a_1\ova_1a_2\ova_2$ and 
$\VV=a_1\ova_2\ova_2$, $\WW=a_1\ova_2$. Then 
$\bra{\VV,\WW}=0$ but, by \theoc{bigons for a pair} and 
\remc{CL2} every pair of geometric representatives of $\VV$ 
and $\WW$ intersects at least in two points. 
\end{ex}

\start{ex}{list} Let $\OO=a_1a_2\ova_1\ova_2 a_3 a_4 
\ova_3\ova_4$. The following is a list of pairs of cyclic 
words $(\VV,\WW)$ such that $\bra{\VV,\WW}=0$ but each pair 
of representatives of $\VV$ and $\WW$ is not disjoint, in the 
parenthesis there is written  the minimal number of 
intersection points of $\VV$ and $\WW$.  (We thank Vladimir Chernov
for pointing out a typo in our previous version).
\begin{romlist}

\item $\VV=\cc(a_1\ova_2\ova_4a_1\ova_2)$ and $\WW=\cc(a_1\ova_2\ova_4\ova_4),$  $(2)$
\item $\VV=\cc(a_1a_3a_3a_3\ova_2)$ and $\WW=\cc(a_1a_3a_3\ova_2),$  $(4)$
\end{romlist}

\end{ex}
\subsection{A counterexample}\label{Counterexample}

\start{ex}{counter} Let $\OO=a_1a_2\ova_1\ova_2$. So the 
surface  $\S_\OO$associated to $\OO$  is a punctured torus. 
Let $\WW=\cc(a_1a_1 a_2 a_2)$. Then $\cob(\WW)=0$. On the 
other hand, by \theoc{bigons} every geometric representative 
of $\WW$ has at least one self-intersection point. More 
generally, for every pair of pair of integers $i,j$, if  
$\WW=\cc(a_1^{i}a_2^j)$ then every representative of $\WW$ 
has at least $(i-1)(j-1)$ self-intersection points and 
$\cob(\WW)=0$. 
\end{ex}

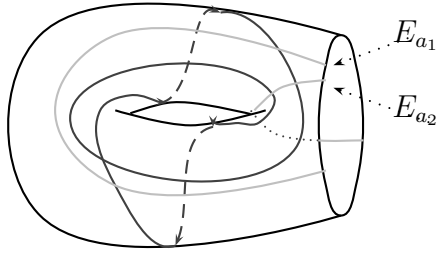
\begin{figure}[ht]
\begin{pspicture}(15.5,4)
\TorusWithOneHole \CurvesForTorusWithOneHole \EdgesForTorusWithOneHole
\rput(10,3){$E_{a_1}$}\rput(10,2){$E_{a_2}$}

\psline[linestyle=dotted, linecolor=black]{->}(10,3)(8.9,2.6)

\psline[linestyle=dotted, linecolor=black]{->}(10.28,1.8)(8.9,2.3)
\end{pspicture}\caption{A representative of $\cc(a_1a_1 a_2 a_2)$ in the punctured
torus}\label{counterexample}
\end{figure}

\start{ex}{more}Let $\OO=a_1a_2\ova_1\ova_2 a_3 a_4 
\ova_3\ova_4$. Then $S_\OO$ is a punctured surface of genus 
two. Here is a sample of cyclic words with cobracket zero, 
with the minimal number of self-intersection points of the 
representatives written in  parenthesis:   
$\cc(a_3a_4\ova_3a_4)$,($1$); 
$\cc(a_2a_3a_2a_3\ova_1\ova_1\ova_1)$, ($2$); 
$\cc(a_3a_1\ova_2a_3a_1\ova_2a_3a_1\ova_2\ova_2\ova_2)$, 
$(2)$; 
 $\cc(\ova_2\ova_2 a_1 a_1 a_1 a_1a_1\ova_2 a_1 a_1 a_1 
a_1a_1)$, $(8)$; $\cc(\ova_2  a_3 a_4a_4 
a_4a_4a_4a_1\ova_2a_3a_1 )$,$(4)$; $\cc(a_3a_1a_3a_1a_2a_2)$ 
($1$).

\end{ex}

Using Examples~\ref{counter} and \refc{more}, we can 
construct examples of classes with cobracket zero which are 
not powers of simple classes, in any surface of genus at 
least one (with or without boundary) except the (closed) 
torus. Observe that, every free homotopy class on the torus 
is a multiple of a simple class, and so the cobracket is 
identically zero.

\subsection{Open problems}
We could not find  examples of non-simple free homotopy 
classes in the pair of pants with cobracket zero. This leaves 
open the following.

\begin{proclama}{Question.} Let $\Sigma$ be two sphere with
three or more punctures. Is every class on $\Sigma$ with cobracket zero  a multiple of a
simple class?
\end{proclama}

When we told Turaev about the above examples, he asked the following question:

\begin{proclama}{Question.}  Are there further ``secondary" operations necessary
to detect simple loops?
\end{proclama}

Our computations suggest an answer for Turaev's question:

\begin{proclama}{A conjectural characterization of simple curves.}
A primitive reduced cyclic word $\VV$ has a simple representative if and only if
$\bra{\VV,\ov{\VV}}=0$.\end{proclama}

By running our program, we could see that for a surface of genus two and one puncture, as
well as for the pair of pants, for all the reduced cyclic words $\VV$ of length at most
fifteen, the bracket of $\VV$ with $\ov{\VV}$ has exactly twice as many terms as the
minimal number of self-intersection points of $\VV$. 
Hence we have the following.

\start{theo}{new}(1) On the sphere with three punctures all 
the cyclic words with at most sixteen letters, except the 
multiples of three peripheral curves, have non-zero 
cobracket. 

(2) On the torus with two punctures all the cyclic words 
$\alpha$ with at most fifteen letters have the property  that 
twice the minimal number of self-intersection points equals 
the number of terms of the bracket 
$[\alpha,\overline{\alpha}]$ in the natural basis. 

\end{theo} 

\begin{proclama}{Question.} Let $\VV$ be a primitive reduced cyclic word.
Is the number of terms of the bracket of $\VV$ with $\ov{\VV}$ equal to twice the minimal
number of self-intersection points of $\VV$?
\end{proclama}

Some examples we have computed suggest that even a more general result may hold:

\begin{proclama}{Question.} Let $n$ and $m$ be two different non-zero integers and let
$\VV$ be a primitive reduced cyclic word. Is the number of terms of the bracket of
$\VV^n$ with $\VV^m$ equal to $2|m.n|$ multiplied by the minimal number of
self-intersection points of $\VV$?\end{proclama}

\appendix
\section{Definition of Lie bialgebra}\label{definition of the Lie bialgebra}

Let $V$ denote a vector space. In order to recall the definition of Lie bialgebra we need
 two auxiliary linear maps
 $\Smap{\om}{V\otimes V\otimes V}$ and $\Smap{\sss}{V\otimes V}$ defined
by $\om(\UU \otimes \VV\otimes \WW) =\WW \otimes \UU \otimes \VV$ and $\sss(\VV\otimes
\WW)=\WW\otimes \VV$ for each triple of cyclic words $\UU,\VV, \WW \in V$.

A {\em Lie algebra} on a vector space $V$ is given by a linear map
$\map{\bri{\,,}}{V\otimes V}{V}$ such that $\bri{\,,} \circ \sss=-\bri{\,,}$ ({\em skew
symmetry}) and $\bri{\,,}(\id \otimes \bri{\,,})(\id+\om+\om^2)=0$ ({\em Jacobi
identity}). A {\em Lie coalgebra} on $V$ is given by a linear map $\map{\cib}{V}{V\otimes
V}$ such that $\sss \circ\cib=-\cib$ ({\em co-skew symmetry}) and $(\id+\om+\om^2)(\id
\otimes \cib)\cib=0$ ({\em cojacobi identity}). $(V,\bri{\,,},\cib)$ is a {\em Lie
bialgebra} if $(V,\bri{\,,})$ is a Lie algebra, $(V,\cib)$ is a Lie coalgebra and the
compatibility equation $\cib \bri{a,b}=\bri{\cib{a},b}+\bri{a,\cib{b}}$ holds for every
$a, b\in V$ where $\bri{a,b\otimes c}=-\bri{b\otimes c,a}=\bri{a,b}\otimes
c+b\otimes\bri{a,c}$. $(V,\bri{\,,},\cib)$ is an {\em involutive Lie bialgebra} if
$(V,\bri{\,,},\cib)$ is a Lie bialgebra and $\bri{\,,}\circ\cib=0$ on $\V$

\begin{rem} According to \cite{Et}, the equation $\bri{\,,}\circ\cib=0$ is the  infinitesimal analogue of having an antipodal
map on a Hopf algebra with square equal to the identity.
\end{rem}

\section{The Lie bialgebra of curves on a surface}\label{curves}

Let $\W$ denote the vector space generated by all (trivial and non-trivial) free homotopy
classes of loops on a surface and let $\W_0$ denote the vector subspace generated by the
class of the trivial loop. In \cite{Gol}, Goldman defined a Lie algebra structure on
$\W$. Turaev proved that this structure passes to the quotient $\W/\W_0$ and that there
exists a cobracket on $\W/\W_0$ which makes the whole structure a Lie bialgebra on
$\W/\W_0$. Notice that the quotient vector space $\W/\W_0$ is canonically isomorphic to
$\V$. If $w \in \W$, we denote by $\conj{w}_0$ the equivalence class of $\W/\W_0$
containing $w$.

We now recall the constructions of Goldman and Turaev. Let 
$\VV$ and $\WW$ be two non-trivial free homotopy classes of 
curves. Choose representatives, $A$ of $\VV$ and $B$ of $\WW$ 
in general position. Hence the intersection of $A$ and $B$ 
consists in a finite number of double points, $p_1,p_2,\dots, 
p_m$. To each of such points  $p_i$ we assign a free homotopy 
class and a sign: the class that contains the loop that 
starts at $p_i$, runs around $A$, and then around $B$, and 
the sign  obtained by comparing the orientation of the 
surface with the orientation given by the branches of $A$ and 
$B$ (in that order) coming out of $p_i$ (positive if it 
coincides, negative otherwise). The bracket $\bri{\VV, \WW}$ 
is defined to be the sum of all the signed free homotopy 
classes.  In symbols, $$\bri{\VV, \WW}=\sum_{p\in A\cap B} 
\sign_p(A,B) \conj{\class(A\cdot_p B)}_0$$ where 
$\sign_p(A,B)$ denotes the sign, $A\cdot_p B$ means the usual 
multiplication of based loops at $p$ and $\class(X)$ denotes 
the conjugacy class of a loop $X$. 

In order to define the cobracket, for each non-trivial free homotopy class $\WW$, choose
a representative $A$ in general position with respect to itself. Thus, its
self-intersection points are finitely many double points $q_1, q_2, \dots, q_m$. To each
of these points $q_i$ one associates an ordered pair of free homotopy classes
$(\WW_1^{q_i},\WW_2^{q_i})$ as follows. Firstly, order the two branches of $A$ coming out
of $q_i$ in such a way that they define the same orientation as the surface.
$\WW_1^{q_i}$ is the conjugacy class of the loop that starts on the first branch, and
runs along $A$ until it arrives to $q_i$ again. Analogously, $\WW_2^{q_i}$ is the
conjugacy class of the loop that starts at the second branch and runs along $A$ until it
finds $q_i$ again.

The cobracket of $\WW$, is given by the formula: $$\cib(\WW)=\sum_{q_1,q_2,
\dots,q_m}\{\WW_1^{q_i}\}_0\otimes \{\WW_2^{q_i}\}_0-\{\WW_2^{q_i}\}_0\otimes
[\WW_1^{q_i}\}_0 .$$

\start{prop}{involutive} The Goldman-Turaev Lie bialgebra is involutive, i.e.,
$\bri{\,,}\circ \cib=0$.
\end{prop}
\begin{proof}
 Let $\WW$ be a free homotopy class and let $\alpha$ be a 
geometric representative of $\WW$ in general position. 
$\bri{\,,}\, \circ \cib(\WW)$ is a sum of free homotopy 
classes over certain ordered pairs $(p,q)$ of 
self-intersection points of $\alpha$. Denote this set of 
pairs of self-intersection points  by $D_\alpha$. Thus, 
$(p,q) \in D_\alpha$ if and only if the arcs of the circle 
between the two preimages of $p$ contain exactly one preimage 
of $q$. Hence for each $(p,q) \in D_\alpha$, there are four 
arcs, $\alpha_1, \alpha_2, \alpha_3$ and $\alpha_4$ of 
$\alpha$ such that $\alpha_1$ and $\alpha_3$ go from $p$ to 
$q$, $\alpha_2$ and $\alpha_4$ go from $q$ to $p$ and 
$\alpha$ runs along $\alpha_1$ then  $\alpha_2$, $\alpha_3$ 
and finally, along $\alpha_4$. 

Clearly, $(p,q) \in D_\alpha$ if and only if $(q,p)\in D_\alpha$. We claim that for each
$(p,q) \in D_\alpha$, the two terms of $\bri{\,,} \circ \cib(\WW)$ corresponding to
$(p,q)$ cancel with the two terms corresponding to $(q,p)$.

The two terms of $\cib(\WW)$ corresponding to $p$ are:
$s_{1,3}(\class(\alpha_1\alpha_2)\otimes
\class(\alpha_3\alpha_4)-\class(\alpha_3\alpha_4)\otimes 
\class(\alpha_1\alpha_2)),$ where  $s_{1,3}=1$ if the 
orientation given by the tangent vector of  $\alpha_1$ at $p$ 
and  the tangent vector of $\alpha_3$ at $p$ coincides with 
the orientation of the surface and $s_{1,3}=-1$ otherwise. 

Now, the terms of bracket of the above linear combination corresponding to $q$ are
$2s_{1,3}s_{2,4} \class(\alpha_4\alpha_3\alpha_2\alpha_1) $ where $s_{2,4}=1$ if the
orientation given by the tangent vector of $\alpha_2$ at $q$ and the tangent vector of
$\alpha_4$ at $q$ coincides with the orientation of the surface and $s_{2,4}=-1$
otherwise.

The same computation for $(q,p)$ gives $-2s_{1,3}s_{2,4}
\class(\alpha_4\alpha_3\alpha_2\alpha_1),$ so the claim is 
proved and so also the proposition. 

\end{proof}

\enddocument